\documentclass[11pt, twoside, letterpaper]{article}
\usepackage{tikz}
\usepackage{verbatim}
\usepackage{color}
\usepackage{amsmath}
\usepackage{amssymb}
\usepackage{amsthm}
\usepackage{bbm}
\usepackage{algpseudocode}

\newtheorem{theorem}{Theorem}[section]
\newtheorem{lemma}{Lemma}[section]

\newtheorem{corollary}{Corollary}[section]

\numberwithin{equation}{section}

\usepackage[margin=1in]{geometry}

\title{Explicit solutions of the SI and Bass models on sparse Erd\H{o}s-R\'enyi and regular networks}
\author{Gadi Fibich, Yonatan Warman}

\date{\today}

\begin{document}

	\date{}
	\maketitle

	\vspace{-3em} 
	\begin{center}
		\footnotesize
		Department of Applied Mathematics, Tel Aviv University \\ fibich@tau.ac.il, warmany98@gmail.com
	\end{center}
	
	\begin{abstract}
		\footnotesize
		We derive explicit expressions for the expected adoption and infection level in the Bass and SI models, respectively, on sparse Erd\H{o}s-R\'enyi networks and on $d$-regular networks. These expressions
	are soloutions of first-order ordinary differential equations, which are fairly easy to analyze. To prove that these expressions are exact, we show that the effect of cycles vanishes as the network size goes to infinity.
	\end{abstract}

\section{Introduction}

Mathematical models of stochastic spreading on networks have found widespread applications across diverse fields of research. In epidemiology, these models simulate the propagation of infectious diseases within social networks~\cite{Epidemics-on-Network-17}. In the social sciences, they provide a framework for analyzing the diffusion of innovations, ideas, and information within social networks~\cite{Barabasi-16,Jackson-08}.
  While stochastic spreading models have been widely applied, the  relationship between the network structure and the dynamics of spreading processes is not well understood. 
  To address this gap, researchers have sought to derive explicit analytical solutions for the spreading dynamics on different types of network structures, as such solutions can provide valuable insights into the key network features that drive or hinder the spreading.

   In this study, we focus on two related network models: the Bass model and the Susceptible-Infected (SI) model. The Bass model describes the adoption of new products or innovations within a population. In this framework, all individuals start as non-adopters and can transition to becoming adopters due to two types of influences: External influences by mass media, and internal influences by individuals who have already adopted the product
   on their peers. 
       The SI model is used to study the spreading of infectious diseases within a population. In this model, some individuals are initially infected (the "patient zero" cases), all subsequent infections occur through internal influences, whereby infected individuals transmit the disease to their susceptible peers, and infected individuals remain contagious indefinitely.
 
   The final outcome of the Bass and SI models is straightforward ---  
   everyone adopts in the Bass model and everyone that can be reached from a patient zero will become infected. The key question is how fast does it take to get to the final outcome, or, more generally, what is 
   the expected adoption/infection level in the population as a function of time.
   Exact explicit expressions for the expected adoption/infection level were obtained for 
   dense networks (infinite homogeneous complete networks, infinite heterogeneous complete networks  that consist of $K$~homogeneous groups, percolated infinite complete networks) and for infinite sparse networks without cycles 
   (homogeneous infinite lines, percolated infinite lines)~\cite{OR-10,Bass-monotone-convergence-23,DCDS-23,Bass-percolation-20,Niu-02}. 
   The absence of cycles on sparse networks leads to a considerable 
   simplification in the analysis. Indeed, Sharky et al.~\cite{Sharkey-15} 
   demonstrated for the SIR model that when the network has a tree structure, the closure of the master (Kolmogorov) equations at the level of triplets is exact. 
   Kiss et al.~\cite{Kiss-15} extended this approach to some special networks with cycles, such as star-triangle networks and networks with cycles of maximum size three. 
   %We note, however, that the closure of the master equations at the level of triplets can lead to systems of ODES that are exact, 
   %but may still be too complex to provide analytical insights easily.

      In this paper, we obtain exact explicit solutions for the Bass and SI models on two types of infinite random networks: Sparse Erd\H{o}s-R\'enyi networks 
      and $d$-regular networks. Both networks are sparse, and contain cycles of any possible length, 
      such that the number of cycles of length~$L$ increases exponentially with~$L$. 
      Our starting point are the master equations.
      We analyze these equations using two analytic tools: 
      The {\em indifference principle} and the {\em funnel theorem}.
      The indifference principle identifies edges that can be added
      or deleted to the network without influencing the nonadoption probability of a given set of nodes, thereby allowing to compute a desired probability on a simpler network~\cite{Bass-boundary-18}. The funnel theorem relates the nonadoption probability of a node of degree~$d$ to the product of the $d$~nonadoption probabilities in which the node is only influenced by one of its edges~\cite{Funnel-25}. This relation is an equality 
      when a node does not lie on any cycle, and a strict inequality when it does.
      
         The explicit expression for the expected adoption/infection level 
          on sparse Erd\H{o}s-R\'enyi networks 
         and $d$-regular networks
         can be computed under the assumption that one can neglect the effect of cycles, i.e., that the funnel relation is an equality for all nodes. We rigorously justify this assumption,
         by utilizing the upper bound of the funnel theorem, which shows that the effect of a cycle on the adoption probability of a node decays as~$1/L^L$, where~$L$ is the cycle length. Therefore, although the number of cycles of length~$L$ increases exponentially with~$L$, the super-exponential decay rate implies that the overall effect of cycles  on the expected adoption/infection level goes to zero as the network size becomes infinite.
     
       On sparse Erd\H{o}s-R\'enyi networks, 
       the analysis of the master equations leads to a single ODE for the 
       expected adoption/infection level. We stress that this ODE is exact, and not an approximation. This ODE is fairly easy to analyze,
       and can be used to show e.g., that the expected adoption/infection level is monotonically decreasing in the average degree, even when the weight of the edges is inversely proportional to the average degree.  
        This monotonicity implies that on any sparse Erd\H{o}s-R\'enyi network,
        the expected adoption/infection level is bounded from below by that of a network of isolated nodes, and from above by a complete infinite network.  
        In the case of the Bass model on sparse Erd\H{o}s-R\'enyi networks, there is no phase transition, and all 
        the nodes ultimately become adopters. The situation is different in the SI model on sparse Erd\H{o}s-R\'enyi networks, where the final infection level as $t \to \infty$ is strictly below one.   
         We derive an algebraic equation for the final infection level, and use it to show that there is no phase transition if the fraction of ``patients zero'' is held fixed as
       the population size becomes infinite.  A phase transition, emerges, however, if the number of ``patients zero'' is held fixed as the population size becomes infinite. 
         
        The analysis of the master equations for $d$-regular networks uses a similar methodology, and also leads to a single ODE for the expected adoption/infection level. This ODE is also fairly easy to analyze. For example, one can use it to show that the expected adoption/infection level increases with~$d$, even when the weights of the edges scales as~$1/d$. 
        We also compare the expected adoption/infection level on $2D$-regular networks and on $D$-dimensional Cartesian networks, 
        as this comparison reveals the difference in the spreading between a deterministic structure and a completely random one.

            From a methodological perspective, the main contribution of this manuscript is the rigorous reduction of the infinite system of master equations into a single ODE. This derivation has two components: Justification that the effect of cycles vanishes as the network size become infinite, and the reduction of the infinite system of master equations into a single ODE. In this manuscript, we applied this methodology to the Bass and SI models on 
            Erd\H{o}s-R\'enyi networks and to $d$-regular networks. 
            We believe that  this approach  has the potential to be extended to other types of sparse networks with cycles (e.g., scale-free networks), to other models (e.g., SIS, SIR, Bass-SIR~\cite{Epidemics-on-Network-17,Bass-SIR-model-16}),
            to the Bass and SI models on hypernetworks~\cite{Hyper-24}, etc.

\section{Bass/SI model on networks}

We begin by introducing a single network model that combines the Bass model for the spreading of new products and the SI model for the spreading of epidemics.
Consider a network~$\cal N$ with $M$~individuals, denoted by ${\cal M}:=\{1, \dots, M\}$. Let $X_j(t)$ denote the state of individual~$j$ at time~$t$, so that 
\begin{equation*}
	X_j(t)=\begin{cases}
		1, \qquad {\rm if}\ j\ {\rm is \ adopter/infected \ at\ time}\ t,\\
		0, \qquad {\rm otherwise,}
	\end{cases}
	\qquad j \in \cal M.
\end{equation*} 
The initial conditions at $t=0$ are stochastic,\footnote{{\em Deterministic initial conditions} are a special case where
	$I_j^0 \in \{0,1\}$.} so that 
\begin{subequations}
	\label{eqs:Bass-SI-models-ME}
	\begin{equation}
		\label{eq:general_initial}
		X_j(0)=	X_j^0 \in \{0,1\}, \qquad j\in {\cal M},
	\end{equation}
	where
	\begin{equation}
	\label{eq:Bass-ER-IC}
		\mathbb{P}(X_j^0=1) =I_j^0, \quad 
		\mathbb{P}(X_j^0=0) =1-I_j^0, 
		\quad I_j^0 \in [0, 1],
		\qquad 
		j \in \cal M,
	\end{equation}
	and
	\begin{equation}
		\label{eq:p:initial_cond_uncor-two_sided_line}
		\mbox{the random variables $\{X_j^0 \}_{j \in \cal M}$ are independent}.
	\end{equation}
So long that $j$ is susceptible, its adoption/infection rate at time~$t$ is
	\begin{equation}
		\label{eq:lambda_j(t)-Bass-model-heterogeneous-tools}
		\lambda_j(t) = p_j+\sum\limits_{k\in {\cal M}} q_{k,j} X_{k}(t),
		\qquad j \in {\cal M}.
	\end{equation}
	Here, $p_j$ is the rate of external influences by mass media on~$j$, and~$q_{k,j}$ is the rate of 
internal influences by~$k$ on~$j$, provided that $k$ is already an adopter/infected. Once~$j$ becomes an adopter/infected, it remains so at all later times.\,\footnote{i.e., the only admissible transition is 
	$X_j=0 ~\longrightarrow ~ X_j=1$.}
Hence, as $ \Delta t \to 0$,
	\begin{equation}
		\label{eq:general_model}
		\mathbb{P} (X_j(t+\Delta  t )=1  \mid   {\bf X}(t))=
		\begin{cases}
			\lambda_j(t) \, \Delta t , &  {\rm if}\ X_j(t)=0,
			\\
			1,\hfill & {\rm if}\ X_j(t)=1,
		\end{cases}
		\qquad 	j \in {\cal M},
	\end{equation}
	where ${\bf X}(t) := \{X_j(t)\}_{j \in \cal M}$ 
	is the state of the network at time~$t$,
	and
	\begin{equation}
		\label{eq:Bass-SI-models-ME-independent}
		\mbox{the random variables $\{X_j(t+\Delta  t )  \mid   {\bf X}(t) \}_{j \in \cal M}$ are independent}.
	\end{equation}  
	%
	%
	% the conditional random variables $\{X_j(t+\Delta  t )=1  \mid   {\bf X}(t) \}_{j \in \cal M}$ are independent.
\end{subequations}
In the Bass model, there are no adopters when the product is first introduced into the market, and so $I_j^0 \equiv 0$.
In the SI model, there are only internal influences for $t>0$, and so $p_j \equiv 0$.

The quantity of most interest is the expected adoption/infection level
at time $t$ in network~${\cal N}$
\begin{equation}
	\label{eq:number_to_fraction-general}
	f(t;{\cal N}):= 
%	\frac{1}{M} \mathbb{E}[N(t)]
%	=
%	\frac{1}{M} \sum_{j=1}^{M}\mathbb{E}[X_j(t)] =
	\frac{1}{M} \sum_{j=1}^{M} f_j(t;{\cal N}),
\end{equation}
where 
%$N:=\sum_{k\in {\cal M}} X_k$ is the number of adopters/infected in the network, and 
$f_j :=\mathbb{P}(X_j = 1) = \mathbb{E}[X_j]$ in the adoption/infection probability of~$j$.
%In the Bass model there are no adopters at $t=0$ when the new product is introduced into the market, and so $I_j^0 \equiv0$. In the SI model, 
%all the infections for $t>0$ 
%are due to interactions with people who were already infected, and so 
%$p_j \equiv 0$.  

\section{Bass/SI model on Erd\H{o}s-Rényi networks}
\label{sec:ER-networks}

In this section, we analyze the Bass/SI model on infinite sparse ER networks.

\subsection{Sparse Erd\H{o}s-Rényi networks}

Let ${\cal G}^{\rm ER}(M,\frac{\lambda}{M})$ denote the set of all undirected graphs with
$M$~nodes, such that  
for any two nodes $ k,j\in{\cal M}$, the edge between~$k$ and~$j$ exists with probability~$\frac{\lambda}{M}$, independently of all other edges, where $0<\lambda <\infty$.   The average degree in~${\cal G}^{\rm ER}(M,\frac{\lambda}{M})$ is 
 	\begin{equation}
 	\label{eq:E_G(degree)}
 	\mathbb{E}_{\cal G} [ \text{degree}(j)] 	=\frac{\lambda (M-1)}{M},
 	\qquad j \in \cal M.
 \end{equation}

Consider a graph $G \in {\cal G}^{\rm ER}(M,\frac{\lambda}{M})$, let ${\bf E} = (e_{k,j})$ denote its adjacency matrix, let all the nodes have weight~$p$, and all the edges have weight~$\frac{q}{\lambda}$.\footnote{
The normalization by~$\lambda$ is chosen so that the sum of the edges weights 
remains  unchanged as $\lambda$ varies, i.e., 
	%$\mathbb{E}_{\cal G}\left[\frac1M \sum_{k,j=1}^M  q_{k,j} \right] $
	\begin{equation}
		\label{eq:E_G{q_j}=q*lambda}
		\begin{aligned}
			\mathbb{E}_{\cal G}\Big[\sum_{k=1}^M \sum_{j=1}^M q_{k,j} \Big] 
			& = \frac{q}{\lambda}  \, \mathbb{E}_{\cal G}\Big[\sum_{k=1}^M \sum_{j=1}^M e_{k,j} \Big] 
			= \frac{q}{\lambda}  \, \mathbb{E}_{\cal G}\Big[ \sum_{k=1}^M \text{\rm deg}(k) \Big] 
			%\\ &  
			=  \frac{q}{\lambda} \lambda(M-1)  \equiv  (M-1) q.
		\end{aligned}
	\end{equation}
		}
The corresponding ER network, denoted by~${\cal N}^{\rm ER}(G)$, is 
\begin{subequations}
	\label{eqs:ER-network}
	\begin{equation}
	\label{eq:ER-network}
	p_j \equiv p, \qquad 
	q_{k,j}=\frac{q}{\lambda}\,e_{k,j},\qquad k,j\in{\cal M}.
	%\label{eq:first}	
\end{equation} 
In addition, let all the nodes have the  initial condition
\begin{equation}
	\label{eq:Bass-model-homog-complete-IC}
	I_j^0 \equiv I^0, \qquad j \in \cal M.
\end{equation}
We  assume that the parameters satisfy 
%\begin{subequations}
%	\label{eqs:assumptions-Bass/{\rm SI}}
\begin{equation}
	0 \le I^0<1, \quad 	q>0, \quad  p \ge 0,
\end{equation}
and 
\begin{equation}
	I^0>0  \qquad  \text{or} \qquad  p>0. 
\end{equation}
\end{subequations}
Specifically, in the Bass model $p>0$ and~$I^0=0$, and 
in the SI model $p=0$ and~$0<I^0<1$.

%The Bass and SI models  on~${\cal N}^{\rm ER}$ are given by, see~\eqref{eqs:Bass-stoch}, 
%\begin{subequations}
%	\label{eqs:Bass-ER}
%	\begin{equation}
%		%	\label{eq:Bass-ER-in}
%		\mathbb{P}\left(X_{j}(0)=1\right)=I^0,\quad
%		\mathbb{P}\left(X_{j}(0)=0\right)=1-I^0,\qquad j\in{\cal M},
%		\label{eq:Bass-ER-IC}
%	\end{equation}
%	the adoption rate of node~$j$ is
%	\begin{equation}
%		\lambda_{j}^{{\rm ER}}(t)=p+\frac{q}{\lambda}\sum_{k \in \cal M}e_{k,j}X_{k}(t),
%	\end{equation}
%	and  as $\Delta t\rightarrow0$,
%	\begin{equation}
%		\mathbb{P}\left(X_{j}(t+\Delta t)=1\left|\,{\bf X}(t)\right.\right)=
%		\begin{cases}
%			\lambda_{j}^{{\rm ER}}(t) \, \Delta t, & \ {\rm if} \  X_{j}(t)=0,\\
%			1, & \ {\rm if} \ X_{j}(t)=1,
%		\end{cases}
%		\qquad \quad  j\in{\cal M}.
%		\label{eq:stam-1-2}
%	\end{equation}
%\end{subequations}
The expected adoption/infection level in the network~${\cal N}^{\rm ER}$ is~$f(t;{\cal N}^{\rm ER})$, see~\eqref{eq:number_to_fraction-general}.
The average of~$f$ over all ER networks in~${\cal G}^{\rm ER}(M,\frac{\lambda}{M})$ is 

	\begin{equation}
	\label{eq:E_G-ER}
	% \langle h \rangle:=
	\mathbb{E}_{{\cal G}}[f](t) :=\sum_{G \in{\cal G}^{\rm ER}(M,\frac{\lambda}{M})} \mathbb{P}^{\rm ER}(G)\, f(t;{\cal N}^{\rm ER}(G)),
	\quad 	\mathbb{P}^{\rm ER}(G)=\left(\frac{\lambda}{M}\right)^{L(G)}\left(1-\frac{\lambda}{M}\right)^{{M \choose 2}-L(G)},
\end{equation}
where $\mathbb{P}^{\rm ER}(G)$
%\begin{equation}
%	\mathbb{P}^{\rm ER}(G)=\left(\frac{\lambda}{M}\right)^{L(G)}\left(1-\frac{\lambda}{M}\right)^{{M \choose 2}-L(G)},
%	\label{eq:E_G-ER}
%\end{equation}
is the probability of the graph~$G$ in~${\cal G}^{\rm ER}(M,\frac{\lambda}{M})$,
and~$L(G)$ is the number of edges in~$G$.

%  The nonadoption level at time~$t$ for a specific
%realization of~\eqref{eqs:Bass-ER} 
%is $\frac{1}{M}\sum_{m=1}^{M}\mathbbm{1}_{X_{m}(t)=0}.$
%Therefore, 

\subsection{Theoretical results}

The first main result of this manuscript is the explicit calculation of the expected adoption/infection level 
 $
f^{\rm ER}:=\lim_{M \to \infty} \mathbb{E}_{\cal G} [f]
$
on infinite sparse ER networks.
We stress that {\em this result is exact}, and does not involve any approximation.
\begin{theorem}
	\label{thm:f^ER}
	The expected adoption/infection level in the Bass/{\rm SI} model~{\rm (\ref{eqs:Bass-SI-models-ME},\ref{eqs:ER-network})} on sparse infinite  {\rm ER}~networks
	is given by
	\begin{subequations}
		\label{eqs:f^ER}
		\begin{equation}
			%	\lim_{M \to \infty}	\overline{\left[S\right]}(t)
			f^{\rm ER}(t;p,q,\lambda,I^0) = 1-
			(1-I^0)e^{-pt-\lambda(1-y(t))}, \qquad t \ge 0,
			\label{eq:f^ER}
		\end{equation}
		where $y(t)$ is the solution of the equation
		\begin{equation}
			\frac{dy}{dt} =\frac{q}{\lambda}\Big(-y+(1-I^0)e^{-pt-\lambda(1-y)}\Big), ~~ t \ge 0,
			\qquad  y(0)=1.
			\label{eq:y-ER}
		\end{equation}
	\end{subequations}
\end{theorem}
 \begin{proof}
	See Section~\ref{sec:proof-thm-ER}.
 \end{proof}

As~$\lambda$ increases, the average degree increases, see~\eqref{eq:E_G(degree)}, but the weight of the edges decreases, see~\eqref{eq:ER-network}, so that the sum of the edges weights remains 
unchanged, see~\eqref{eq:E_G{q_j}=q*lambda}. 
Therefore, it is not obvious whether $f^{\rm ER}(t;\lambda)$ should  
increase or decrease with~$\lambda$.
The explicit expression~\eqref{eqs:f^ER} enables us
to prove that $f^{\rm ER}$ increases with~$\lambda$:
\begin{lemma}
	\label{lem:f^ER-increasing-in-lambda}
	Let $t>0$. Then $f^{\rm ER}(t;\lambda)$ is monotonically increasing in~$\lambda$.
\end{lemma}
 \begin{proof}
See Appendix~\ref{app:f^ER-increasing-in-lambda}.
\end{proof}

We can ulilize the monotonicity of~$f^{\rm ER}$ in~$\lambda$ to obtain tight lower and upper bounds for the expected adoption/infection level in all sparse ER networks:
\begin{lemma}
	\label{lem:f^ER_Bass-bounds}
	Consider the Bass/{\rm SI} model~{\rm (\ref{eqs:Bass-SI-models-ME},\ref{eqs:ER-network})} on infinite sparse {\rm ER}~networks. 
	\begin{enumerate}
		\item 
		When $\lambda=0$, the expected/infection adoption level is the same as on networks of isolated nodes,
		i.e., 
		\begin{equation}
			\label{eq:f^ER_lambda=0}
			f^{\rm ER}(t;p,q,I^0,\lambda=0) = f^{\rm isolated}(t;p,I^0),
			\qquad 	f^{\rm isolated}:=I^0 + (1-I^0)(1-e^{-pt}).
		\end{equation}

		\item  %As $\lambda \to \infty$, 
	%	the expected adoption/infection level approaches 
		%that of the  compartmental Bass/{\rm SI} model~\eqref{eq:df_dt-compartmental}	i.e., 
%		\begin{equation}
%			\label{eq:lim_lambda_to_infty_f^ER}
		$\lim_{\lambda \to \infty} f^{ \text{\rm ER}}(t;p,q,I^0,\lambda) = f^{\rm compart}(t;p,q, I^0)$,
		%\end{equation}
		where $ f^{\rm compart}$ is the expected adoption/infection level
		in the  compartmental Bass/{\rm SI} model
		\begin{equation}
			\label{eq:f^complete_infty-ODE}
			\frac{df}{dt} = (1-f)(p+qf), \qquad f(0) = I^0.
		\end{equation}

%\begin{equation}
%	\label{eq:f^compart}
%	f^{\rm compart}:=
%	\begin{cases}
%		\vspace{3mm}
%		\displaystyle\frac{\frac{1+\frac{q}{p} I^0}{1-I^0}-e^{-(p+q)t}}{\frac{1+\frac{q}{p} I^0}{1-I^0}+\frac{q}{p}e^{-(p+q)t}}, & \ \text{\rm if} \ p>0, \\
%		\displaystyle\frac{1}{1+\big(\frac1{I^0}-1\big)e^{-qt}} , & \ \text{\rm if} \ p=0,
%	\end{cases}
%$\end{equation}	
	\item  For $0<\lambda<\infty$,
		\begin{equation}
			\label{eq:bounds-f^ER}
			f^{\rm isolated}(t;p,I^0)  < f^{\rm ER}(t;p,q,I^0,\lambda)  < f^{\rm compart}(t;p,q,I^0), \qquad t>0.
		\end{equation}
	\end{enumerate}
\end{lemma}
%\begin{lemma}
%	\label{lem:f^ER_Bass-bounds}
%	Consider the Bass/{\rm SI} model~{\rm (\ref{eqs:Bass-SI-models-ME},\ref{eqs:ER-network})}  on infinite sparse  {\rm ER}~networks. Then
%	$$
%I^0 + (1-I^0)(1-e^{-pt})  < f^{\rm ER}(t;p,q,I^0,\lambda)  < f^{\rm compart}(t;p,q,I^0), \qquad t>0, \quad  0<\lambda<\infty,
%	$$
%	where 
%\begin{equation}
%\label{eq:f^compart}
%f^{\rm compart}:=
%\begin{cases}
%	\vspace{3mm}
%	\displaystyle\frac{\frac{1+\frac{q}{p} I^0}{1-I^0}-e^{-(p+q)t}}{\frac{1+\frac{q}{p} I^0}{1-I^0}+\frac{q}{p}e^{-(p+q)t}}, & \ \text{\rm if} \ p>0, \\
%	 \displaystyle\frac{1}{1+\big(\frac1{I^0}-1\big)e^{-qt}} , & \ \text{\rm if} \ p=0,
%\end{cases}
%\end{equation}
%%and $A:=\frac{1+\frac{q}{p} I^0}{1-I^0}$.
%\end{lemma}
\begin{proof} See Appendix~\ref{app:f^ER_Bass-bounds}. 
\end{proof}
 
%This result is intuitive. The lower bound $f^{\rm ER}(\cdot,\lambda=0)=I^0 + (1-I^0)(1-e^{-pt})$ is the expected adoption/infection level 
%in a network of isolated modes. The upper bound~$\lim_{\lambda \to \infty} f^{\rm ER}(\cdot,\lambda) = f^{\rm compart}$ is  the expected adoption/infection level in
%infinite complete homogeneous networks, see~\cite{Bass-monotone-convergence-23}. 

The $\lambda  \to \infty$ limit is intuitive,  
since $f^{\rm compart}$ is the limit as $M \to \infty$ of the expected adoption/infection level on complete homogeneous networks
as the degree of the nodes goes to infinity and the weight of the edges goes to zero~\cite{Bass-monotone-convergence-23,Niu-02}.
This is also the case for sparse ER~networks, where the weight of the edges $\frac{q}{\lambda}$ goes to zero
and the average degree goes to infinity  as~$\lambda \to \infty$.

\subsubsection{Bass model on sparse Erd\H{o}s-Rényi networks}

The Bass model on sparse ER~networks is given by~(\ref{eqs:Bass-SI-models-ME},\ref{eqs:ER-network}) with  $p>0$ and $I^0 = 0$. 
By Theorem~\ref{thm:f^ER},  the expected adoption level 
$f^{\rm ER}_{\rm Bass}(t):=f^{\rm ER}(t;I^0 = 0)$  is given by
\begin{subequations}
	\label{eqs:f^ER-Bass}
	\begin{equation}
		\label{eq:f^ER-Bass}
		%	\lim_{M \to \infty}	\overline{\left[S\right]}(t)
		f^{\rm ER}_{\rm Bass}(t) = 1-
		e^{-pt-\lambda(1-y_{\rm Bass}(t))}, \qquad t \ge 0,
	\end{equation}
	where $y_{\rm Bass}(t):=y(t;I^0 = 0)$ is the solution of the equation
	\begin{equation}
		\frac{dy}{dt} = - \frac{q}{\lambda}\left(y-e^{-pt-\lambda(1-y)}\right), \quad t > 0,
		\qquad  y(0)=1.
		\label{eq:y-ER-Bass}
	\end{equation}
\end{subequations}

%In the Bass model on sparse ER networks, $p>0$ and $I^0 = 0$.  
%%In that case, the explicit expression for~$f^{\rm ER}$ has the following representation: 
%\begin{corollary}
%	Let $f^{\rm ER}_{\rm Bass}(t;p,q,\lambda):=f^{\rm ER}(t;p,q,\lambda,I^0 = 0)$
%	denote the expected adoption level in the Bass model~{\rm (\ref{eqs:Bass-SI-models-ME},\ref{eqs:ER-network})} on infinite sparse ER networks. Then
%	\begin{subequations}
%		\label{eq:f^ER-Bass}
%		\begin{equation}
%			%\lim_{M\rightarrow\infty}\overline{\left[S\right]}=
%			f^{\rm ER}_{\rm Bass}(t;p,q,\lambda) = 1-
%			e^{-pt-\lambda(1-y_{\rm Bass}(t))}, \qquad t \ge 0,
%			\label{eq:final-bass}
%		\end{equation}
%		where $y_{\rm Bass}(t)$ is the solution of the equation
%		\begin{equation}
%			\label{eq:y-ER-Bass}
%			\frac{dy}{dt}=\frac{q}{\lambda}\left(-y+e^{-pt-\lambda(1-y)}\right), ~~ t \ge 0,
%			\qquad  y(0)=1.
%		\end{equation}
%	\end{subequations}
%\end{corollary}
% \begin{proof}
%	This follows from Theorem~\ref{thm:f^ER} by substituting $I^0=0$. 
% \end{proof}

By Lemma~\ref{lem:f^ER-increasing-in-lambda}, $f^{\rm ER}_{\rm Bass}(t)$ is monotonically increasing in~$\lambda$, and by Lemma~\ref{lem:f^ER_Bass-bounds}, 
$$
	1-e^{-pt}  < f^{\rm ER}_{\rm Bass}(t)  <\frac{1-e^{-(p+q)t}}{1+\frac{q}{p}e^{-(p+q)t}}, 
	\qquad t>0, \quad  0<\lambda<\infty.
	$$
	Note that the upper bound~$\frac{1-e^{-(p+q)t}}{1+\frac{q}{p}e^{-(p+q)t}}$ 
	is the expected adoption level in the compartmental Bass model~\cite{Bass-69}.
	
Using the explicit solution~\eqref{eqs:f^ER-Bass}, we can show that $f^{\rm ER}_{\rm Bass}(t)$ increases monotonically with time, and that all the network nodes become adopters as $t \to \infty$: 
\begin{lemma}
	\label{lem:f^ER_Bass-monotone}
	The expected adoption level~$f^{\rm ER}_{\rm Bass}(t)$ is monotonically increasing in~$t$ from $f^{\rm ER}_{\rm Bass}(0)=0$ to
	$f^{\rm ER}_{\rm Bass}(\infty)=1$.
\end{lemma}
\begin{proof}
		See Appendix~\ref{app:f^ER_Bass-monotone}.
\end{proof}
%	See Appendix~\ref{app:f^ER_Bass-monotone}.
% \end{proof}

%We can ulilize the monotonicity in~$\lambda$ to obtain tight lower and upper bounds for the expected adoption level in the Bass model on all sparse ER networks:
%\begin{lemma}
%	\label{lem:f^ER_Bass-bounds}
%	Consider the Bass model on infinite sparse ER networks. Then
%	$$
%	1-e^{-pt}  < f^{\rm ER}_{\rm Bass}(t;p,q,\lambda)  < f_{\rm Bass}(t;p,q), \qquad t>0, \quad  0<\lambda<\infty,
%	$$
%	where $f_{\rm Bass}:=\frac{1-e^{-(p+q)t}}{1+\frac{q}{p}e^{-(p+q)t}}$ is the expected adoption level in the Bass model on infinite complete networks, see~\eqref{eq:f_Bass-Niu}.
%\end{lemma}
% \begin{proof} See Section~\ref{app:f^ER_Bass-bounds}. 
% 

%\subsection{Simulations}

\begin{figure}[ht!]
	\begin{center}
		\scalebox{0.6}{\includegraphics{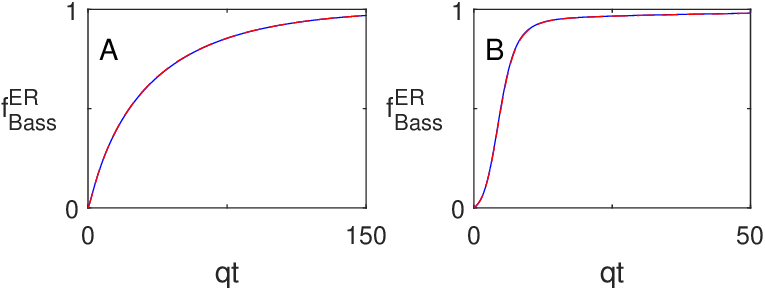}}
	%	\scalebox{0.6}{\includegraphics{f_ER.eps}}
		%		\scalebox{0.5}{\includegraphics{Figure_26.eps}}
	\end{center}
	\caption{The expected adoption level~$f^{\rm ER}_{\rm Bass}$ in the Bass model~{\rm (\ref{eqs:Bass-SI-models-ME},\ref{eqs:ER-network})}
		on sparse infinite ER networks.		
		Dashed orange line is  the numerical
		solution with $M=2000$ nodes, averaged over 100~simulations. 
		Solid blue line is the explicit expression~\eqref{eqs:f^ER-Bass}. The two curves are indistinguishable.
		%	(${\rm A}_2$) The auxiliary function~$y_{\rm Bass}(t)$, see~\eqref{eq:y-ER-Bass}. 
		Here $I^0=0$, $p=0.001$, and $q=0.05$. (A)~$\lambda=\frac{1}{2}$. 
		(B)~$\lambda=3$.}
	\label{fig:f-ER}
\end{figure}

Figure~\ref{fig:f-ER} illustrates the numerical agreement 
between the expected adoption level in simulations of the Bass model on  
sparse ER~networks, and the explicit expression~\eqref{eqs:f^ER-Bass}. Although  $p \ll q$ in these simulations, 
the adoption curve~$f^{\rm ER}_{\rm Bass}(t)$ is concave when $\lambda = \frac12$, which indicates that the adoption process is 
primarily driven by external influences. 
This is because when $\lambda = \frac12$, the fraction of isolated nodes is $\mathbb{P}_{\rm deg}(0)= e^{-\frac12}  \approx 60\%$, and so  most of the adoptions are external.
% and $\mathbb{P}_{\rm deg}(1)= \frac12 e^{-\frac12}  \approx 0.3$,
%see~\eqref{eq:pois}, only 10\% of the nodes have degree $\ge 2$. Consequently, 
%,  and so\,\footnote{%Since  $0<y(t)<1$, see Lemma~\ref{lem:y(t)-properties-ER}, this approximation also follows from~\eqref{eq:final-bass}
	%	By continuity, this approximation also follows from~\eqref{eq:f^ER_lambda=0} for $\lambda \ll 1$. } 
%$$
%		f^{\rm ER}_{\rm Bass}(t)  \approx  1- e^{-pt}.
%$$
When $\lambda=3$, however,  $X(\lambda = 3) \approx 94\%$ of the nodes belong to the giant component. Since, in addition, $q \gg p$, most of the adoptions are internal. Hence, $f^{\rm ER}_{\rm Bass}$~has an S-shape.

%Recall that in sparse ER networks the giant component emerges at $\lambda = 1$ (Section~\ref{sec:ER-sparse}).
%Expression~\eqref{eq:f^ER} for~$f^{\rm ER}$, however, varies smoothly
%in~$\lambda$ . This smoothness  can  be seen in 

\begin{figure}[ht!]
	\begin{center}
		\scalebox{0.5}{\includegraphics{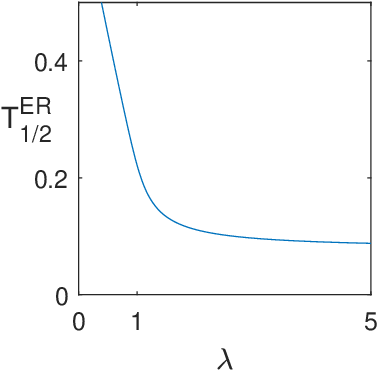}}
		%		\scalebox{0.5}{\includegraphics{Figure_26.eps}}
	\end{center}
	\caption{ The half-life in the Bass model~{\rm (\ref{eqs:Bass-SI-models-ME},\ref{eqs:ER-network})} on sparse infinite ER networks, as a function of~$\lambda$.
		%when the weights of all edges is~$q$ (blue solid)  and $\frac{q}{\lambda}$ (red dashes). 
		Here $I^0=0$, $p=0.001$, and $q = 0.05$.}	
	\label{fig:half_life_lambda}
\end{figure}

%More generally, since the average degree is given by~$\lambda$,  
%see~\eqref{eq:lim_M2infty<deg>=lambda}, and since adding edges accelerates the adoption process 
%(Corollary~\ref{cor:f^A<f^B}), the adoption process becomes 
%faster as $\lambda$ increases. To see that,
In Figure~\ref{fig:half_life_lambda} we plot 
the half-life $T^{\rm ER}_{1/2}(p,q,\lambda)$ in infinite sparse ER networks, 
which is the time where half of the population adopted the product, 
and observe that 
it decreases (i.e., the spreading becomes faster) as~$\lambda$ increases,
in agreement with Lemma~\ref{lem:f^ER-increasing-in-lambda}. 
Note that 
{\em  there is no phase transition at $\lambda=1$}.  This is because 
in the Bass model,  the adoption process depends on 
local properties of the network structure~\cite{OR-10},
and these properties do not undergo a phase transition at $\lambda=1$. 
A phase transition can occur at~$\lambda=1$ 
when the spreading depends on global properties of the network structure,
such as connectivity or average distance between nodes.
Indeed, this can be the case in the SI model on sparse ER networks, see Section~\ref{sec:ER-SI}.

\subsubsection{SI model on sparse Erd\H{o}s-Rényi networks}
\label{sec:ER-SI}

The SI model on sparse ER~networks is given by~(\ref{eqs:Bass-SI-models-ME},\ref{eqs:ER-network}) with $p=0$ and $0<I^0<1$. 
By Theorem~\ref{thm:f^ER},  the expected infection level 
$f^{\rm ER}_{\rm SI}(t):=f^{\rm ER}(t;p = 0)$  is given by
\begin{subequations}
	\label{eqs:f^ER-SI}
	\begin{equation}
		\label{eq:f^ER-SI}
		%	\lim_{M \to \infty}	\overline{\left[S\right]}(t)
		f^{\rm ER}_{\rm SI}(t;q,I^0,\lambda) = 1-(1-I^0)
		e^{-\lambda(1-y_{\rm SI}(t))}, \qquad t \ge 0,
	\end{equation}
	where $y_{\rm SI}(t):=y(t;p = 0)$ is the solution of the equation
	\begin{equation}
		\frac{dy}{dt} = - \frac{q}{\lambda}\left(y-(1-I^0)e^{-\lambda(1-y)}\right), \quad t > 0,
		\qquad  y(0)=1.
		\label{eq:y-ER-SI}
	\end{equation}
\end{subequations}

By Lemma~\ref{lem:f^ER-increasing-in-lambda}, $f^{\rm ER}_{\rm SI}(t)$ is monotonically increasing in~$\lambda$, and by Lemma~\ref{lem:f^ER_Bass-bounds}, 
$$
I^0  < f^{\rm ER}_{\rm SI}(t)  <\frac{1}{1+(\frac1{I^0}-1)e^{-qt}}, 
\qquad t>0, \quad  0<\lambda<\infty.
$$
	Note that the upper bound~$\frac{1}{1+(\frac1{I^0}-1)e^{-qt}}$ 
is the expected infection level in the compartmental SI model.

\begin{figure}[ht!]
\begin{center}
	\scalebox{0.6}{\includegraphics{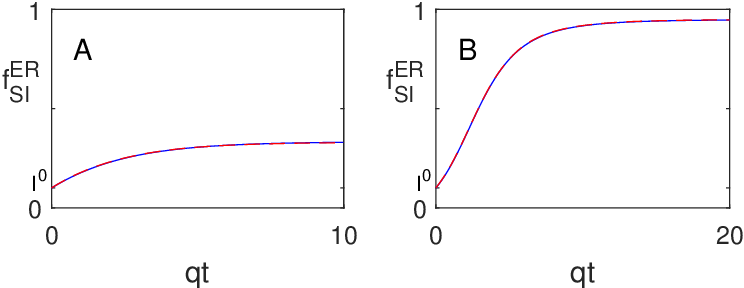}}
	%		\scalebox{0.5}{\includegraphics{Figure_26.eps}}
\end{center}
\caption{The expected infection level~$f^{\rm ER}_{\rm SI}$ in the SI model~{\rm (\ref{eqs:Bass-SI-models-ME},\ref{eqs:ER-network})}
	on sparse infinite ER networks. 
	Dashed orange line is the numerical	solution with $M=2000$~nodes, averaged over 100~simulations.
	Solid blue line is the explicit expression~\eqref{eqs:f^ER-SI}.  
	%		(${\rm A}_2$) The auxiliary function~$y_{\rm SI}(t)$, see~\eqref{eq:y-ER-SI}. 
	%	
	Here $I^0=0.1$ $p=0$, and $q=0.05$. (A)~$\lambda=0.9$. 
	(B)~$\lambda=3$.}	
\label{fig:f_ER_SI}
\end{figure}

Figure~\ref{fig:f_ER_SI}  confirms the numerical agreement between simulations of the SI model on sparse ER networks 
and the explicit expression~\eqref{eqs:f^ER-SI} for~$f^{\rm ER}_{\rm SI}$.  {\em The most notable difference from 
the Bass model is that a positive fraction of the population does not become infected 
as $t\rightarrow\infty$}.
Indeed, we can use the explicit expression~\eqref{eqs:f^ER-SI} 
to show that the limiting infection level~$f_{\rm SI}^{\infty}$ is strictly below one:
\begin{lemma}
\label{lem:f^SI-monotone}
			The expected infection level~$f^{\rm ER}_{\rm SI}(t)$  
in the {\rm SI} model on sparse infinite {\rm ER}~networks
is monotonically increasing in~$t$ from $f^{\rm ER}_{\rm SI}(0) = I^0$ to  $f_{\rm SI}^{\infty}:= 
\lim_{t\rightarrow\infty}f^{\rm ER}_{\rm SI}(t)$,
where  $0<f_{\rm SI}^{\infty}<1$. Furthermore, 
$f_{\rm SI}^{\infty}$~is the unique solution in~$(0, 1)$ of 
\begin{equation}
	\label{eq:implicit} 
	f_{\rm SI}^{\infty} = 
	1- (1-I^0)e^{-\lambda f_{\rm SI}^{\infty}}.
\end{equation}
\end{lemma}
 \begin{proof}
See Appendix~\ref{app:f^ER_SI-monotone}.
 \end{proof}

%Making the change of variables $X:=(1-y_{\rm SI}^{\infty})$, we get the following implicit
%function
%\begin{equation}
%	X= (1-I^0)e^{-\lambda X}.
%	\label{eq:implicit}
%\end{equation}
%This equation can be solved numerically and we get
%\begin{equation}
%	f_{\infty}=1-y_{\infty}=1-(1-I^0)e^{-X}.\label{eq:f_infty}
%\end{equation}

Recall that if the largest connected component of a sparse infinite ER network 
contains a positive fraction of the nodes, it is called the {\em giant component}.
Let~$X = X(\lambda)$ denote the fraction of nodes in the giant component. Then~$X$  
 satisfies the equation
\begin{equation}
\label{eq:GC-ER}
X = 1-e^{-\lambda X}. 
\end{equation}
When $0 \le \lambda \le 1$, the only solution of~\eqref{eq:GC-ER} is $X(\lambda) \equiv 0$,
and so there is no giant component. When $\lambda>1$, however,  
eq.~\eqref{eq:GC-ER} also admits a positive solution, which is the size of the 
giant component. See e.g.,~\cite{Barabasi-16,Networks-10} for further details and proofs.

\begin{figure}[ht!]
\begin{center}
	\scalebox{0.6}{\includegraphics{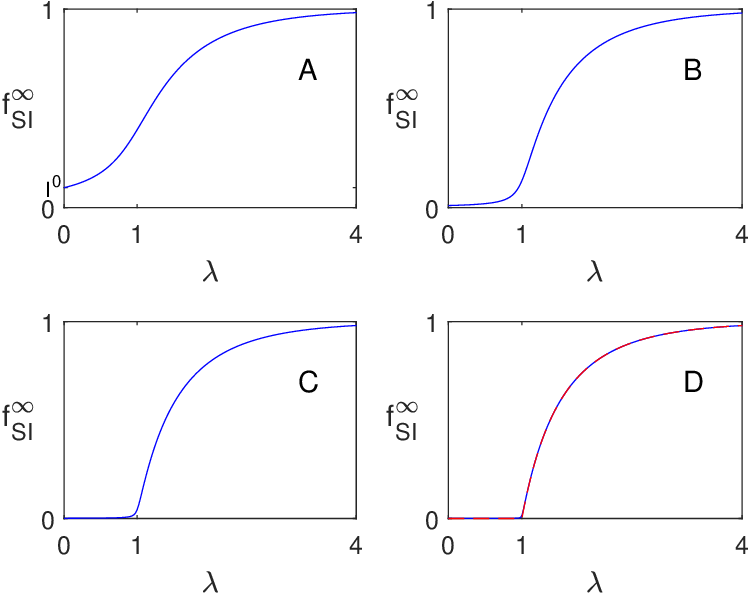}}
	%		\scalebox{0.5}{\includegraphics{Figure_26.eps}}
\end{center}
\caption{Final fraction of infected individuals in the SI model on infinite sparse ER
	networks, as a function of~$\lambda$, see~\eqref{eq:implicit}. 
	The initial fraction of infected is: 
	(A)~$I^0=0.1$, (B)~$I^0=0.01$. (C)~$I^0=0.001$.
	(D)~$I^0=0.0001$. The dashed orange curve is $X(\lambda)$, 
	see~\eqref{eq:GC-ER}. The two curves are indistinguishable.
}
\label{fig:f-ER-infty(t)}
\end{figure}

In Figure~\ref{fig:f-ER-infty(t)} we 
plot the final fraction~$f_{\rm SI}^{\infty}$ of infected nodes
 as a function of~$\lambda$.
When $I^0 = 0.1$, there is no phase transition at $\lambda = 1$.  
As we let $I^0 \to 0+$, however, a phase transition emerges. 
Indeed, we have 
\begin{corollary}
\label{cor:lim_beta_0_to_0_f_SI^infty=X(lambda)}
Let $0  < \lambda<\infty$. Then 
$ \lim_{ I^0 \to 0} f_{\rm SI}^{\infty}(q,I^0,\lambda) = X(\lambda)$, where~$X$ is the solution of~\eqref{eq:GC-ER} which is positive for $\lambda>1$. 
\end{corollary}
 \begin{proof}
 See Appendix~\ref{app:lim_beta_0_to_0_f_SI^infty=X(lambda)}.
\end{proof}

Thus, a phase transition occurs if we hold the number of initially-infected 
nodes constant and let $M \to \infty$, but not if we 
 hold the fraction of initially-infected 
nodes constant and let $M \to \infty$.

Corollary~\ref{cor:lim_beta_0_to_0_f_SI^infty=X(lambda)}  shows that
{\em as $I^0 \to 0$,  the only nodes that become infected are those that belong to the giant component}.
 Indeed, the nodes that are not in the giant component belong to small clusters that have a vanishing probability to include a patient zero as $I^0 \to 0$. 
 Corollary~\ref{cor:lim_beta_0_to_0_f_SI^infty=X(lambda)} 
 is illustrated numerically in Figure~\ref{fig:f-ER-infty(t)}D, where the plot
 of~$f_{\rm SI}^{\infty}(\lambda)$ for $I^0 = 0.0001$ is
 indistinguishable from that of~$X(\lambda)$.

\section{Bass/SI model on $d$-regular  networks}
\label{sec:d-regular-networks}

In this section, we analyze the Bass/SI model on $d$-regular networks.

\subsection{$d$-regular networks}

    A graph is called ``$d$-regular'', if the degree of all its nodes is equal to~$d$. 
	The family~${\cal G}_{M,d}^{\rm regular}$ of random $d$-regular graphs
	consists of all $d$-regular  graphs with $M$~nodes, such 
	that the probability of any graph in~${\cal G}_{M,d}^{\rm regular}$ is the same.
	Consider a graph $G \in {\cal G}_{M,d}^{\rm regular}$, let ${\bf E} = (e_{k,j})$ denote its adjacency matrix, let all the nodes have weight~$p$, and all the edges have weight~$\frac{q}{d}$.\,\footnote{
		The normalization by~$d$ is chosen so that the sum of the edges weights of each node 
		remains  unchanged as $d$~varies, i.e., 
		%$\mathbb{E}_{\cal G}\left[\frac1M \sum_{k,j=1}^M  q_{k,j} \right] $
		\begin{equation}
			\label{eq:q_j=q-d-regular}
		\sum_{k=1}^M q_{k,j} 
				 = \frac{q}{d} \sum_{k=1}^M e_{k,j} = q , \qquad j \in \cal M.
		\end{equation}
	}
	The corresponding $d$-regular network, denoted by~${\cal N}^{ \text{\rm d-reg}}={\cal N}^{ \text{\rm d-reg}}(G)$, is 
	\begin{subequations}
		\label{eqs:d-regular-network}
		\begin{equation}
			\label{eq:d-regular-network}
			p_j \equiv p, \qquad 
			q_{k,j}=\frac{q}{d}\,e_{k,j},\qquad k,j\in{\cal M}.
			%\label{eq:first}	
		\end{equation} 
		In addition, let all the nodes have the  initial condition
		\begin{equation}
			\label{eq:Bass-model-homog-complete-IC-d-reguler}
			I_j^0 \equiv I^0, \qquad j \in \cal M.
		\end{equation}
		We  assume that the parameters satisfy 
		%\begin{subequations}
		%	\label{eqs:assumptions-Bass/{\rm SI}}
		\begin{equation}
			0 \le I^0<1, \qquad 	q>0, \quad  p \ge 0,
		\end{equation}
		and 
		\begin{equation}
			I^0>0  \qquad  \text{or} \qquad  p>0. 
		\end{equation}
	\end{subequations}

The expected adoption/infection level in network~${\cal N}^{ \text{\rm d-reg}}$ is
$
f(t;{\cal N}^{ \text{\rm d-reg}}) := \frac{1}{M}\sum_{j=1}^{M} f_j(t;{\cal N}^{ \text{\rm d-reg}})$,
and the average of~$f$ over all $d$-regular networks  
is $\mathbb{E}_{{\cal G}_{M,d}^{\rm regular}} [f]$. 

\subsection{Theoretical results}

The second  main result of this manuscript is the explicit calculation of the expected adoption/infection level 
%$
%f^{ \text{\rm d-reg}}(t):= \lim_{M \to \infty}\mathbb{E}_{{\cal G}_{M,d}^{\rm regular}} [f]
%$
on infinite $d$-regular networks. This result is also exact, and does not involve any approximation.
\begin{theorem}
	\label{thm:f^d-regular}
		Let 
	\begin{equation}
		\label{eq:f^d-regular-def}
		f^{ \text{\rm d-reg}}(t):= \lim_{M \to \infty}\mathbb{E}_{{\cal G}_{M}^{\text{\rm d-reg}}} [f]
	\end{equation}
	denote the expected adoption/infection level in the Bass/{\rm SI }model~{\rm (\ref{eqs:Bass-SI-models-ME},\ref{eqs:d-regular-network})} on $\lim_{M \to \infty} {\cal G}_{M}^{\text{\rm d-reg}}$.
	Then 
%	
%		The expected adoption/infection level in the Bass/{\rm SI} model~{\rm (\ref{eqs:Bass-SI-models-ME},\ref{eqs:d-regular-network})} on  infinite $d$-regular networks
%	is given by
	\begin{subequations}
		\label{eqs:f^d-regular}
		\begin{equation}
			\label{eq:f^d-regular=1-[S^d-reg]}
			f^{ \text{\rm d-reg}} = 1- [S^{ \text{\rm d-reg}}],
		\end{equation}
		where $[S^{ \text{\rm d-reg}}]$ is the solution of the equation
		\begin{equation}
			\label{eq:S^d-reg}
			\frac{d [S]}{dt} = -[S] \left( p+q \Big(1- 
			%	e^{-\frac{2}{d}pt}  (1-I^0)^{\frac{2}{d}} [S]^{\frac{d-2}{d}} 
			\Big(\frac{[S]}{e^{-pt}(1-I^0)}\Big)^{-\frac{2}{d}} [S]
			\Big) \right),
			\qquad [S](0 )= 1-I^0.
		\end{equation}
	\end{subequations}
\end{theorem}
\begin{proof}
	See Section~\ref{sec:f^d-regular-proof}.
\end{proof}

The availability of the explicit expression~\eqref{eqs:f^d-regular} enables us
to show that $f^{ \text{\rm d-reg}}(t)$ is monotonically increasing in~$d$.
This result is not obvious, since as $d$ increases each node 
can be influenced by more nodes, but the influence rate of each of these nodes becomes weaker,  so that the sum of the edges weights of each node remains 
unchanged, see~\eqref{eq:q_j=q-d-regular}. 
\begin{lemma}
	\label{lem:f^d-regular-increasing-in-d}
	Let $t>0$. Then $f^{ \text{\rm d-reg}}(t)$ is monotonically increasing in~$d$.
\end{lemma}
\begin{proof} See Appendix~\ref{app:f^d-regular-increasing-in-d}.
\end{proof}

Lemma~\ref{lem:f^d-regular-increasing-in-d}  is of the same spirit as Lemma~\ref{lem:f^ER-increasing-in-lambda}. 
Together, they show that, roughly speaking, {\em few strong links lead to a lower adoption/infection level that numerous strong ones}.

We can ulilize the monotonicity of~$f^{ \text{\rm d-reg}}$ in~$d$ to obtain tight lower and upper bounds
for the expected adoption/infection level in all $d$-regular networks:
\begin{lemma}
	\label{lem:bounds-f^regular}
	Consider the Bass/{\rm SI} model~{\rm (\ref{eqs:Bass-SI-models-ME},\ref{eqs:d-regular-network})} on  infinite $d$-regular networks. 
	\begin{enumerate}
		\item 
		When $d=2$, the expected adoption level is the same as on infinite circles/lines,
		i.e., 
		\begin{equation}
			\label{eq:f^2-regular_equiv_f^1D}
			f^{\rm 2-regular}(t;p,q,I^0) \equiv f^{\rm 1D}(t;p,q,I^0),
		\end{equation}
		%   where $f^{\rm 1D}$ is the expected adoption level in the Bass model on  infinite circles,
		%  see~\eqref{eqs:f^1D}. 
		where
				\begin{equation}
			\label{eq:f^1D}
			f^{\rm 1D}:=
			\begin{cases}
				1-(1-I^0) e^{-(p+q)t+ q(1-I^0)\frac{1-e^{-pt}}{p}}, &\ \text{\rm if} \ p>0,\\
				1-(1-I^0) e^{-qI^0 t
					
				}, & \ \text{\rm if} \ p=0.
			\end{cases}
		\end{equation}

			\item  
%		As $d \to \infty$, the expected adoption level approaches that on  infinite complete networks,
%		 i.e., 
%		\begin{equation}
%			\label{eq:lim_d_to_infty_f^d-regular}
			$\lim_{d \to \infty} f^{ \text{\rm d-reg}}(t;p,q,I^0) = f^{\rm compart}(t;p,q, I^0)$,
			where $ f^{\rm compart}$ is the expected adoption/infection level
		in the  compartmental Bass/{\rm SI} model~\eqref{eq:f^complete_infty-ODE}.
		\item  For $d=2,3, \dots$,
		\begin{equation}
			\label{eq:bounds-f^regular}
			f^{\rm 1D}(t;p,q,I^0) \le f^{ \text{\rm d-reg}}(t;p,q,I^0)  < f^{\rm compart}(t;p,q, I^0), \qquad t>0.
		\end{equation}
	\end{enumerate}
\end{lemma}

\begin{proof} See Appendix~\ref{app:bounds-f^regular}.
\end{proof}

The result for 2-regular networks is to be expected, since 2-regular networks are collections of disjoint cycles, whose average size goes to infinity as $M \to \infty$. The $d \to \infty$ limit is also intuitive,  
since $f^{\rm compart}$ is the limit as $M \to \infty$ of the expected adoption/infection level on complete networks
as the degree of the nodes goes to infinity and the weight of the edges goes to zero~\cite{Bass-monotone-convergence-23,Niu-02}.

\begin{figure}[ht!]
	\begin{center}
		\scalebox{0.6}{\includegraphics{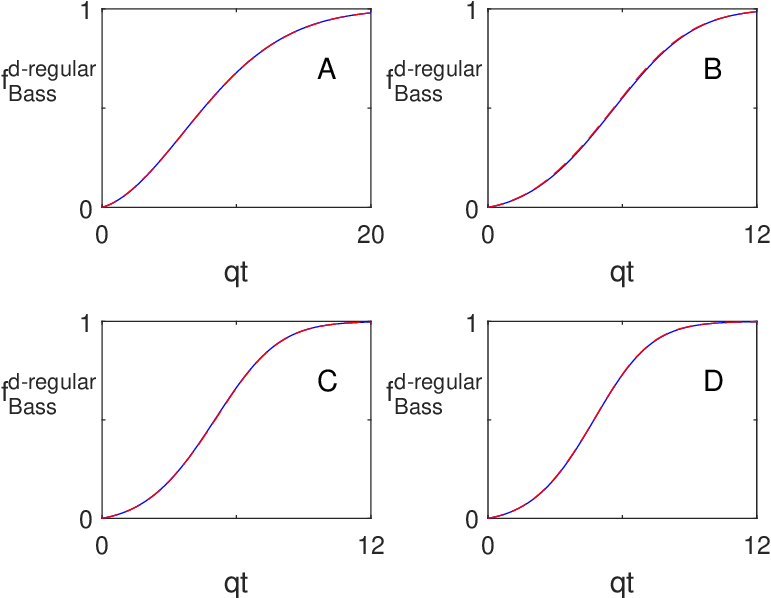}}
		%		\scalebox{0.5}{\includegraphics{Figures/ER/Figure_26.eps}}
	\end{center}
	\caption{The expected adoption level in the Bass model on $d$-regular networks.		
		The dashed line is the numerical
		solution of~{\rm (\ref{eqs:Bass-SI-models-ME},\ref{eqs:d-regular-network})} with $M=10,000$ nodes, averaged over 10~simulations. 
		The solid line is the explicit expression~\eqref{eqs:f^d-regular}. Here $p=0.001$, $q=0.05$, and $I^0=0$.  A)~$d=2$. 
		B)~$d=3$.  C)~$d=4$.  D)~$d=5$. 
	}	
	\label{fig:f-regular}
\end{figure}

\begin{figure}[ht!]
	\begin{center}
		\scalebox{0.6}{\includegraphics{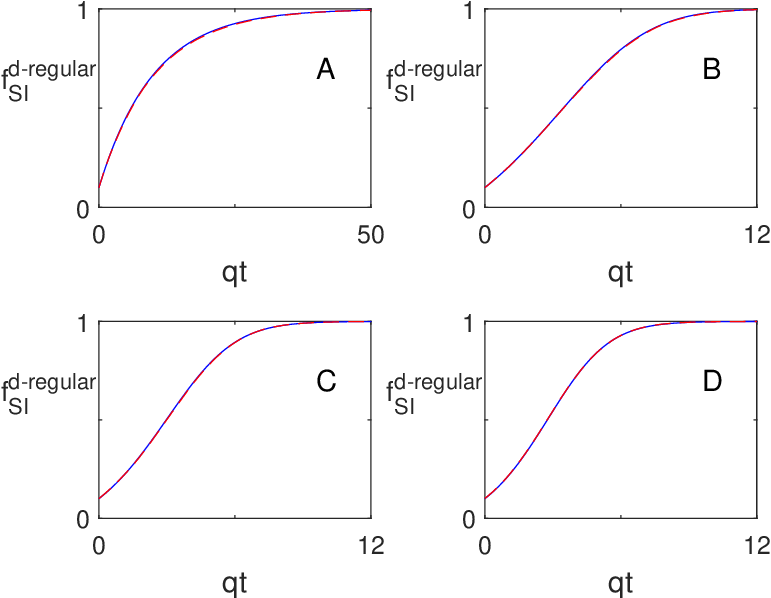}}
	\end{center}
	\caption{Same as Figure~\ref{fig:f-regular} for the SI model.
	Here $p=0$, $q=0.05$, and $I^0=0.1$.
	}	
	\label{fig:f-regular-SI}
\end{figure}

\begin{figure}[ht!]
	\begin{center}
		\scalebox{0.7}{\includegraphics{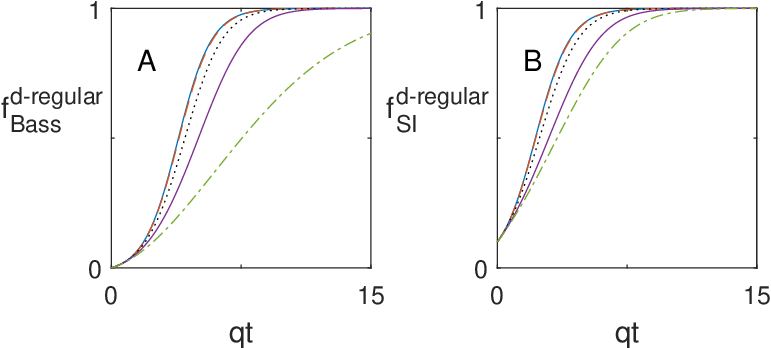}}
	\end{center}
	\caption{The expected adoption/infection level $f^{ \text{\rm d-reg}}(t)$ in the Bass/SI model~{\rm (\ref{eqs:Bass-SI-models-ME},\ref{eqs:d-regular-network})}
		 on $d$-regular networks, for $d=2$ 
		(dash-dot), $d=4$ (solid), $d=10$ (dots), and $d=100$ (dashes).
		The $d=100$ curve is nearly indistinguishable from~$f^{\rm compart.}$ (solid), 
		see~\eqref{eq:f^complete_infty-ODE}.  A)~Bass model with $p=0.001$, $I^0 = 0$, and $q=0.05$. 
		B)~SI model with $p=0$, $q=0.05$, and $I^0 = 0.1$.
	}	
	\label{fig:f-regular-compare}
\end{figure}

Figures~\ref{fig:f-regular} and~\ref{fig:f-regular-SI} illustrate the numerical agreement between the explicit expression~\eqref{eqs:f^d-regular} for~$f^{ \text{\rm d-reg}}(t)$ and simulations  of the Bass and SI models, respectively, on $d$-regular networks.
The $d$-regular graphs in these simulations were generated using the algorithm of Steger and Wormald~\cite{Steger-99}. 
Figure~\ref{fig:f-regular-compare} 
%presents $f^{ \text{\rm d-reg}}(t)$ for several values of~$d$, and 
confirms the predictions of Lemma~\ref{lem:f^d-regular-increasing-in-d} that 
$f^{ \text{\rm d-reg}}(t)$ is monotonically increasing in~$d$, 
and that it approaches~$f^{\rm compart}(t)$ as 
$d \to \infty$.

\subsection{Comparison of $d$-regular and Cartesian networks}

Let~$f^{\rm D}(t)$ denote the expected adoption/infection level in the Bass/SI model on infinite $D$-dimensional Cartesian networks~$\mathbb{Z}^D$, where each node is connected to its 
$2D$ nearest neighbors and the weight of all the edges is~$\frac{q}{2D}$. 
It is instructive to compare~$f^{\rm D}(t)$ with~$f^{ \text{\rm 2D-reg}}(t)$, 
since in both networks all the nodes have the same degree and all the edges have the same weight.
Therefore, this comparison reveals {\em  the difference in the spreading 
	between a deterministic Cartesian structure and a 
	completely random one}.
	Note that this difference also provides an upper bound for the effect of adding a {\em small-worlds structure}~\cite{Watts-98} to a Cartesian network.

\begin{figure}[ht!]
	% created using cellular_1D_2D_3D_4D_k_regular_main.m
	\begin{center}
		\scalebox{0.6}{\includegraphics{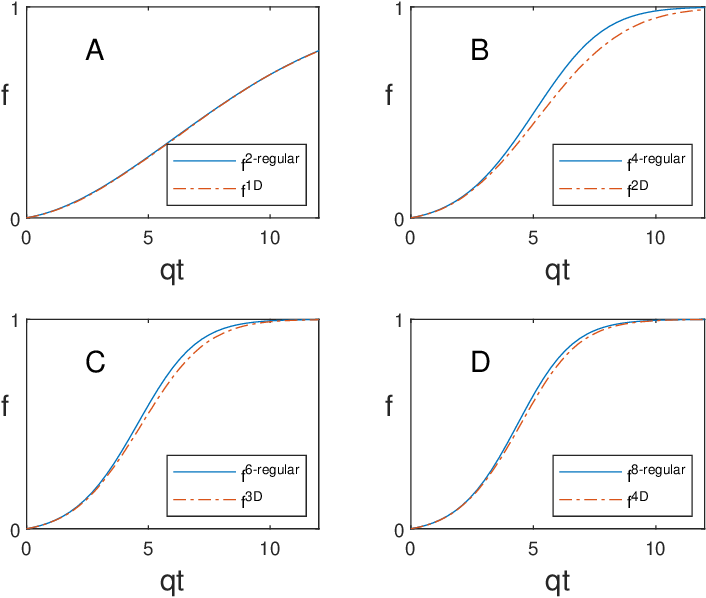}}
	\end{center}
	\caption{The expected adoptions levels $f^{\rm D}_{\rm Bass}$ (dash-dot) and $f^{ \text{\rm 2D-reg}}_{\rm Bass}$ (solid) in the Bass model on  $D$-dimensional Cartesian networks and on $2D$-regular networks.
		Here $p = 0.002$, $q=0.1$, and $I^0=0$. 
		A)~$D=1$.  
		B)~$D=2$.  
		C)~$D=3$.  
		D)~$D=4$.  
	}	
	\label{fig:cellular_1D_2D_3D_4D_k_regular}
\end{figure}

\begin{figure}[ht!]
	% created using cellular_1D_2D_3D_4D_k_regular_main.m
	\begin{center}
		\scalebox{0.6}{\includegraphics{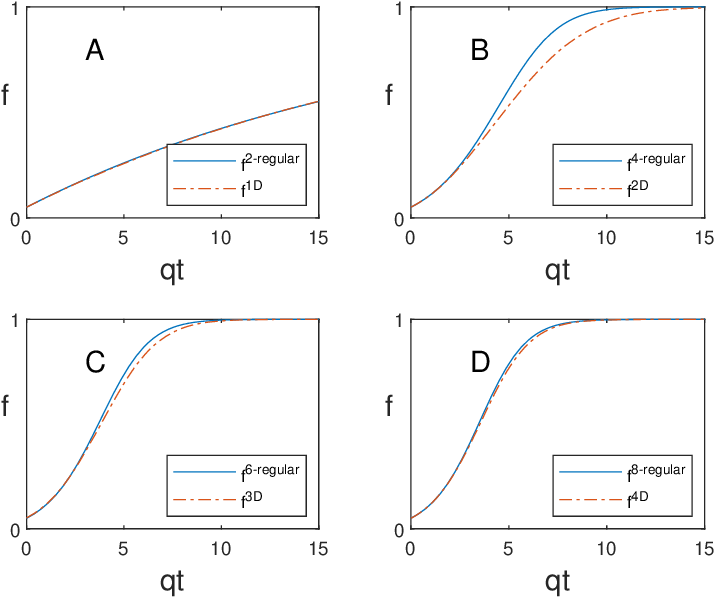}}
	\end{center}
	\caption{Same as Figure~\ref{fig:cellular_1D_2D_3D_4D_k_regular} for the SI model. 
		Here $p = 0$, $q=0.1$, and $I^0=0.05$. 
	}	
	\label{fig:cellular_1D_2D_3D_4D_k_regular_SI}
\end{figure}

The two adoption levels are identical when $D=1$ and 
as $D \to \infty$, see Lemma~\ref{lem:bounds-f^regular}, i.e., 
\begin{equation}
\label{eq:f^d-regular-f^D}
f^{ \text{\rm 2-regular}}(t) = f^{\rm 1D}(t) ,
\qquad \lim_{d \to \infty} f^{ \text{\rm d-reg}}(t)  =  f^{\rm compart}(t)  = \lim_{D \to \infty} f^{\rm D}(t) .
\end{equation}
Numerical simulations (Figures~\ref{fig:cellular_1D_2D_3D_4D_k_regular}
and~\ref{fig:cellular_1D_2D_3D_4D_k_regular_SI})
suggest that for $1<D<\infty$,
$$
f^{ \text{\rm 2D-reg}}(t)> f^{\rm D}(t), \qquad t>0.
$$
This result is intuitive. Indeed, the creation of external adopters/infected on both networks is identical. Any external adopter/infected evolves into a cluster of adopters/infected through internal influences. 
On Cartesian networks, the cluster roughly expands as a $D$-dimensional sphere~\cite{OR-10}. This is not the case on $d$-regular networks, since there is no spatial correlation between the cluster nodes (beyond the edges along which it expands). As a result, the cluster nodes in a Cartesian network are less likely to be in direct contact with susceptible nodes than the cluster nodes in a $d$-regular network. Hence, the expansion rate of clusters on  $d$-regular networks
is faster than on $D$-dimensional networks.

Figures~\ref{fig:cellular_1D_2D_3D_4D_k_regular}
and~\ref{fig:cellular_1D_2D_3D_4D_k_regular_SI} also show that 
the difference between~$f^{\rm D}$  and~$f^{ \text{\rm 2D-reg}}$ decays as 
$D \to 1$  and as  $D \to \infty$, in accordance 
with~\eqref{eq:f^d-regular-f^D}. 
The difference between~$f^{\rm D}$  and~$f^{ \text{\rm 2D-reg}}$ increases  
as $\frac{p}{q}$ and/or $I^0$ decrease (data not shown), since under these parameters changes the internal influences become more dominant than the external ones, 
and so the network structure has a larger effect of the spreading.

\section{Explicit calculation  of $f^{\rm ER}$}
\label{sec:proof-thm-ER}

In this section we provide a proof of Theorem~\ref{thm:f^ER}.

\subsection{Preliminaries}

We first recall some results on sparse ER networks, and 
derive an estimate for the number of cycles of length~$L$. %This estimate will be used in Lemma~\ref{lem:S_B_d-funnel} below.
Let 
$$
d_j:=\text{degree}(j), \qquad j \in \cal M.
$$
The degree distribution in ${\cal G}^{\rm ER}(M,\frac{\lambda}{M})$ is
\begin{equation}
	\mathbb{P}_{\rm deg}^{M}(d):=\mathbb{P}(d_j=d)={M-1 \choose d}\left(\frac{\lambda}M\right)^{d} \left(1-\frac{\lambda}M\right)^{M-1-d}, \qquad d = 0, \dots, M-1.
	\label{eq:binomic}
\end{equation}
%While the degrees of all nodes are 
%identically distributed, they are not independent.  
%For example, if the nodes $\{1, \dots, M-1\}$ have degree zero,  node~$M$ can only have degree zero.
As $M \to \infty$, it approaches 
%The degrees of nodes in $\lim_{M \to \infty}{\cal G}^{\rm ER}(M,\frac{\lambda}{M})$ have 
the Poisson distribution 
\begin{equation}
	\mathbb{P}_{\rm deg}(d):=\lim_{M\rightarrow\infty}\mathbb{P}_{\rm deg}^{M}(d)
	=\frac{\lambda^{d}}{d!}e^{-\lambda},\qquad d=0,1,\dots
	\label{eq:pois}
\end{equation}

Let $j \in \cal M$, and let~$c_L(j;G)$ denote the number of cycles of size~$L$ that go through node~$j$ in a graph~$G$.
% and let
%\begin{equation}
%	B_{d}(G):=\Big\{ m\in{\cal M} \,  \big |  \, 
%	d_m=d   \Big \} 
%	\label{eq:a_k_def}
%\end{equation}
%denote the subset of all the nodes of degree~$d$ in~$G$. 
The expected number of cycles of size~$L$ that go through node~$j$, given that 
its degree is equal to~$d$,  is %denoted by 
%given that $d_j=d$, is 
\begin{equation}
	C_L^d:=
	\mathbb{E}_{{\cal G}(M,\frac{\lambda}{M})}
	\big[c_L(j;G) \mid d_j=d \big],
	%\frac{1}{|B_{d}|}\sum_{m\in B_{d}}c_L(m;G), 
	\qquad d=2,3,\dots
	\label{eq:c_L^d}
\end{equation} 
%where~$|B_{d}|$ is the size of~$B_{d}$. 
Since the family~${\cal G}(M,\frac{\lambda}{M})$ is invariant under permutations of the nodes,
$C_L^d$ is independent of $j$. Therefore, $C_L^d$ 
is also the expected number of cycles of size~$L$ that go through nodes
of degree~$d$. 

\begin{lemma}
	\label{lem:E_G[c_L^d]-bound}
	%in ${\cal N}^{\rm ER} \in {\cal G}(M,\frac{\lambda}{M})$. 
	Let $ L \in \{3, \dots, M\}$. Then  
	%	the expectation of~$c_L^d$	in~${\cal G}(M,\frac{\lambda}{M})$ has the  bounds
	\begin{equation}
		\label{eq:E_G[c_L^d]-bound}
%		\frac{\kappa_{L,d}^{\rm ER}}{M} \left(1-\frac{L-1}{M}\right)^{L-3}
%		<
		C_L^d <
		\frac{\kappa_{L,d}^{\rm ER}}{M}, \qquad \kappa_{L,d}^{\rm ER}:=\frac{ d (d-1)}{2} \lambda^{L-2}.
	\end{equation}
	In addition, if we fix $L$ and let $M \to \infty$, then
	\begin{equation}
		\label{eq:E_G[c_L^d]}
		C_L^d \sim 
		\frac{\kappa_{L,d}^{\rm ER}}{M}, \qquad M \to \infty.
	\end{equation}
\end{lemma}

%Let~$c_L(m;G)$ denote the number of cycles with $L$~nodes in 
%a graph~$G \in {\cal G}^{\rm ER}(M,\frac{\lambda}{M})$ that go through node~$m \in \cal M$, and let
%\begin{equation}
%	B_{d}(G):=\Big\{ m\in{\cal M} ~  \big | ~ 
%	\text{\rm degree}(m)=d   \Big \} 
%	\label{eq:a_k_def}
%\end{equation}
%denote the set of all nodes of degree~$d$. 
%The average number of cycles with $L$~nodes that go through each node of degree~$d$
%is denoted by 
%\begin{equation}
%	c_L^d(G):=\frac{1}{|B_{d}|}\sum_{m\in B_{d}}c_L(m;G), \qquad d=2,3,\dots,
%	\label{eq:c_L^d}
%\end{equation} 
%where~$|B_{d}|$ is the size of~$B_{d}$. 
%
%\begin{lemma}
%	\label{lem:E_G[c_L^d]-bound}
%	%in ${\cal N}^{\rm ER} \in {\cal G}^{\rm ER}(M,\frac{\lambda}{M})$. 
%	For any $ L=3, \dots, M$, 
%	the expectation of~$c_L^d$	has the upper bound
%	\begin{equation}
%		\label{eq:E_G[c_L^d]-bound}
%		\mathbb{E}_{\cal G}[c_L^d] <
%		\frac{C_{L,d}^{\rm ER}}{M}, \qquad C_{L,d}^{\rm ER}:=\frac{ d (d-1)}{2} \lambda^{L-2}.
%	\end{equation}
%	In addition, if we fix $L$ and let $M \to \infty$, we get 
%	\begin{equation}
%		\label{eq:E_G[c_L^d]}
%		\mathbb{E}_{\cal G}[c_L^d] \sim 
%		\frac{C_{L,d}^{\rm ER}}{M}, \qquad M \to \infty.
%	\end{equation}
%\end{lemma}
\begin{proof} See Section~\ref{app:E_G[c_L^d]-bound}.
\end{proof}

\subsection{$\protect\overline{[S^\infty]}=\sum_{d=0}^{\infty}\mathbb{P}_{\rm deg}(d) \, \overline{[S^\infty_{B_{d}}]}$}

The expected susceptibility level in the network~${\cal N}^{\rm ER}={\cal N}^{\rm ER}(G)$ is
\begin{equation*}
[S](t;{\cal N}^{\rm ER}):= 1-f(t;{\cal N}^{\rm ER}) = \frac{1}{M}\sum_{m=1}^{M}[S_m](t;{\cal N}^{\rm ER}),
\qquad 
%\label{eq:check-1-1}
%\end{equation}
[S_m] :=1-f_m.
\end{equation*}

The average of~$[S]$ %and~$[S_{B_{d}}]$
over all the networks in~${\cal G}(M,\frac{\lambda}{M})$ is denoted by, see~\eqref{eq:E_G-ER},
\begin{equation}
	\begin{aligned}
		\overline{[S]}(t):=\mathbb{E}_{{\cal G}(M,\frac{\lambda}{M})}\big [ \, [S](t;{\cal N}^{\rm ER})\big].
		% \qquad
		%	=\sum_{G\in {\cal G}(M,\frac{\lambda}{M})}
		%	 \mathbb{P}^{\rm ER}(G) \, [S]\big(t;{\cal N}^{\rm ER}(G)\big), 
		%	\\
		%	\overline{[S_{B_{d}}]}(t):=\mathbb{E}_{{\cal G}(M,\frac{\lambda}{M})}\big [ \, [S_{B_{d}}](t;{\cal N}^{\rm ER})\big].
		%	=\sum_{G \in {\cal G}(M,\frac{\lambda}{M})}\mathbb{P}^{\rm ER}(G) \, [S_{B_{d}}]\big(t;{\cal N}^{\rm ER}(G)\big).
	\end{aligned}
	\label{eq:E-E-R}	
\end{equation}

%Since the family ${\cal G}(M,\frac{\lambda}{M})$ is invariant under permutations of the nodes,  it follows that 

\begin{lemma}
	For any $j \in \cal M$, 
	\begin{equation}
		\begin{aligned}
			\overline{[S]}(t)=\mathbb{E}_{{\cal G}(M,\frac{\lambda}{M})}\big [ \, [S_j](t;{\cal N}^{\rm ER})\big].
			%	=\sum_{G\in {\cal G}(M,\frac{\lambda}{M})}
			%	 \mathbb{P}^{\rm ER}(G) \, [S]\big(t;{\cal N}^{\rm ER}(G)\big), 
			%	\\
			%	=\sum_{G \in {\cal G}(M,\frac{\lambda}{M})}\mathbb{P}^{\rm ER}(G) \, [S_{B_{d}}]\big(t;{\cal N}^{\rm ER}(G)\big).
		\end{aligned}
		\label{eq:E-E-R-j}	
	\end{equation}
\end{lemma}
\begin{proof}
	The family ${\cal G}(M,\frac{\lambda}{M})$ of sparse ER networks is invariant under permutations of the nodes.
	Therefore, $\mathbb{E}_{{\cal G}}\Big [ [S_m]\Big]$ is independent of~$m$,
	and so 	$
	\overline{[S]} = \mathbb{E}_{{\cal G}}\big [ [S]\big]
	=  \mathbb{E}_{{\cal G}}\Big [ \, \frac1M \sum_{m=1}^M[S_m]\Big]
	=  \frac1M \sum_{m=1}^M  \mathbb{E}_{{\cal G}}\Big [ [S_m]\Big]
	= \mathbb{E}_{{\cal G}}\Big [ [S_j]\Big].
	$
	%	Similarly,
	%$$
	%	\overline{[S_{B_{d}}]}
	%	=\mathbb{E}_{{\cal G}}\big [ \, [S_{B_{d}}]\big]
	%	=\mathbb{E}_{{\cal G}}\Big [ \frac{1}{|B_{d}|}\sum_{m\in B_{d}}[S_m]\Big]
	%	=\mathbb{E}_{{\cal G}}\Big [ \frac{1}{|B_{d}|}\sum_{m=1}^M [S_m]\Big]
	%$$	
	%	
	%	and
	%	$$
	%	\mathbb{E}_{{\cal G}}\big [ \,  [S_{j}] ~ \Big| ~ d_j = d  \big]
	%	$$
\end{proof}

%in~(\ref{eq:stam-1-1}) and $[S_m]$ and in~(\ref{eq:s_by_ak})
%are defined over one fixed ER network~${\cal N}^{\rm ER}$.
%We can 
%
%these probabilities 
%.
Thus, $\overline{[S]}$ is also the expected susceptibility level of a node.
We can similarly define the conditional expected susceptibility level of a node, given that its degree is~$d$, i.e., 
\begin{equation}
	\label{eq:S_Bd_by_dist}
	\overline{[S_{B_{d}}]}:=\mathbb{E}_{{\cal G}(M,\frac{\lambda}{M})}\big [ \,  [S_{j}]~ \Big| ~ d_j = d  \big].	
\end{equation}
Since ${\cal G}(M,\frac{\lambda}{M})$  is invariant under permutations of the nodes,
% and~\eqref{eq:E-E-R-M=infty}, respectively.
$\overline{[S_{B_{d}}]}$ is independent of $j$.

%\subsection{$\protect\overline{[S^\infty]}=\sum_{d=0}^{\infty}\mathbb{P}_{\rm deg}(d) \, \overline{[S^\infty_{B_{d}}]}$}

Let  
\begin{equation}
	\label{eq:E-E-R-M=infty}
	\overline{[S^\infty]}:= \lim_{M \to \infty}\overline{[S]},
	\qquad 	\overline{\left[S^\infty_{B_{d}}\right]}:= \lim_{M \to \infty}\overline{[S_{B_{d}}]}.
\end{equation}
Since
$
f^{\rm ER} = 1-\overline{[S^\infty]},
$
it is enough to compute $\overline{[S^\infty]}$.
To do that, we first prove
\begin{lemma}
	\label{lem:lemma2}
	Consider the Bass/{\rm SI} model~{\rm (\ref{eqs:Bass-SI-models-ME},\ref{eqs:ER-network})} on $\lim_{M\to \infty} {\cal G}(M,\frac{\lambda}{M})$. Then
	\begin{equation}
		\label{eq:s_by_dist}
		\overline{[S^\infty]}(t)=\sum_{d=0}^{\infty}\mathbb{P}_{\rm deg}(d) \,
		%	\lim_{M\rightarrow\infty}
		\overline{[S^\infty_{B_{d}}]}(t),
	\end{equation}
	where $\mathbb{P}_{\rm deg}(d)$ is given by~\eqref{eq:pois}.
	% and\,\footnote{Note that $\overline{[S_{B_{d}}]}$ is independent of $j$, Since ${\cal G}(M,\frac{\lambda}{M})$  is invariant under permutations of the nodes.} 
	% and~\eqref{eq:E-E-R-M=infty}, respectively.
	%		\begin{equation}
		%		\label{eq:S_Bd_by_dist}
		%	\overline{\left[S^\infty_{B_{d}}\right]}:= \lim_{M \to \infty}\overline{[S_{B_{d}}]}, \qquad 
		%		\overline{[S_{B_{d}}]}:=\mathbb{E}_{{\cal G}(M,\frac{\lambda}{M})}\big [ \,  [S_{j}]~ \Big| ~ d_j = d  \big].	
		%	\end{equation}
\end{lemma}
\begin{proof}
	By the law of total expectation, and using~\eqref{eq:E-E-R-j} and~\eqref{eq:S_Bd_by_dist},  
	$$
	\overline{[S]} = \mathbb{E}_{{\cal G}(M,\frac{\lambda}{M})}\big [\, [S_{j}]\, \big]
	=\sum_{d=0}^{M-1}  \mathbb{P}(d_j = d) \,    %\mathbb{P}_{\rm deg}^{M}(d)
	\mathbb{E}_{{\cal G}(M,\frac{\lambda}{M})}\big [ \,  [S_{j}]~ \Big| ~ d_j = d  \big]
	=\sum_{d=0}^{M-1}  \mathbb{P}(d_j = d) \,    %\mathbb{P}_{\rm deg}^{M}(d)
	\overline{[S_{B_{d}}]}.
	$$	 
	Since $ \mathbb{P}(d_j = d)  = \mathbb{P}_{\rm deg}^{M}(d)$, letting $M \to \infty$ and using~\eqref{eq:pois}
	and~\eqref{eq:E-E-R-M=infty} 
	%\eqref{eq:s_by_dist},  and~\eqref{eq:S_Bd_by_dist} 
	yields the result. 
\end{proof}

\subsection{Vanishing effect of cycles}

We recall the {\em Funnel Theorem} for undirected networks: 
\begin{theorem}[\cite{Funnel-25}]
	\label{thm:funnel-ER}
	Consider the Bass/{\rm SI} model~\eqref{eqs:Bass-SI-models-ME}
	on an undirected network~$\cal N$, in which all the nodes have weight~$p$ and initial condition~$I^0$, 
	and all the edges have weight~$q$. 
	Let $j \in \cal M$, 
	such that $d_j:=\text{degree}(j) \ge 2$, 
	and denote by $\{k_{i}\}_{i=1}^{d_j}$ the nodes that are connected to~$j$.
	For $i=1, \dots, d_j$, 
	define the network~${\cal N}^{k_i,p_j}$  by deleting in~${\cal N}$ all the edges of~$j$,
	except for the edge $j \leftrightarrow k_i$,
	and denote $[S^{k_i,p_j}_j]:=[S_j](t;{\cal N}^{k_i,p_j})$.
	Define also the network~${\cal N}^{p_j}$ 
	by deleting in~${\cal N}$ all the edges of~$j$,
	and denote $[S^{p_j}_j]:=[S_j](t;{\cal N}^{p_j})$. 
	\begin{itemize}
		\item If $j$ does not belong to any cycle in $\cal N$, then
		\begin{equation}
			\label{eq:funnel-ineq-ER-d}
			[S_j] =
			\frac{\prod_{i=1}^{d_j} [S^{k_i,p_j}_j]}
			{([S^{p_j}_j])^{d_j-1}},   \qquad \text{\rm \bf (funnel equality)}
		\end{equation} 
		where
		\begin{equation}
			\label{eq:[S^pj_j](t;N^ER)=(1-beta0)e^-pt}
			[S^{p_j}_j]  =  (1-I^0)e^{-pt}.
		\end{equation} 
		\item  If $j$ belongs to the $N_j$~cycles $\{C_n\}_{n=1}^{N_j}$, then 	
		\begin{equation}
			\label{eq:funnel-upper-bound-Ncycle}
			0 < 	[S_{{j}}]  -\frac{\prod_{i=1}^{d_j} [S^{k_i, p_j}_j]}{([S^{p_j}_j])^{d_j-1}}
			< 
			[S_j^{p_j}]	
			\sum_{n=1}^{N_j}
			E(t;L_n), \qquad t > 0  ,	 \qquad \text{\rm \bf (funnel inequality)}
			%		 e^{-\left(p+q\right)t}\left(\frac{eqt}{\left \lfloor{\frac{L_n +1}{2}}\right \rfloor}\right)^{\left \lfloor{\frac{L_n +1}{2}}\right \rfloor},
			%		 \quad    \min_{1 \le n \le N} \left \lfloor{\frac{L_n +1}{2}}\right \rfloor \ge qt,
		\end{equation}
		where $L_n$~is the number of nodes of~$C_n$, and~$E(t;L_n)$ satisfies the bounds
		\begin{subequations}
			\label{eqs:[S_i,j]-[S_i][S_j]-N-paths-global-in-time-E}
			\begin{equation}
				\label{eq:[S_i,j]-[S_i][S_j]-N-paths-global-in-time-E-temporal}
				E(t;L) \le 2(1-I^0)e^{-\left(p+ q \right)t} \bigg(\frac{e q t}{\left \lfloor{\frac{L+1}{2}}\right \rfloor}\bigg)^{\left \lfloor{\frac{L+1}{2}}\right \rfloor}, 
				\qquad    0 \le t  < \frac1q  \bigg \lfloor{\frac{L+1}{2}} \bigg \rfloor, 
				%	\qquad 
				%	E(t;L_n) \le  \left(\frac{q}{p+q}\right \rfloor}\right)^{\left \lfloor{\frac{L_n+1}{2}}\right \rfloor},
			%	\qquad   t \ge 0.
		\end{equation}
		and
		\begin{equation}
			\label{eq:[S_i,j]-[S_i][S_j]-N-paths-global-in-time-E-global}
			E(t;L) \le 2(1-I^0) \left(\frac{q}{p+q} \right)^{\left \lfloor{\frac{L+1}{2}}\right \rfloor},
			\qquad  t \ge 0 . % \quad n=1, \dots, N_j.
		\end{equation}
	\end{subequations}
\end{itemize}
\end{theorem}

The following lemma shows that the effect of cycles on ER networks becomes negligible as $M \to \infty$, so that one can effectively 
compute $\overline{\left[S^\infty_{B_{d}}\right]}$ using the funnel equality:   
\begin{lemma}
\label{lem:S_B_d-funnel}
%Consider the Bass and {\rm SI} models~\eqref{eqs:Bass-ER} on $ \lim_{M \to \infty} {\cal G}^{\rm ER}(M,\frac{\lambda}{M})$. Then
Let $\overline{\left[S^\infty_{B_{d}}\right]}$ be given by~{\em (\ref{eq:S_Bd_by_dist},\ref{eq:E-E-R-M=infty})}. Then 
\begin{equation}
	\overline{\left[S^\infty_{B_{d}}\right]}
	= 		\lim_{M \to \infty}
	\mathbb{E}_{{\cal G}(M,\frac{\lambda}{M})}  \bigg[ \frac{\prod_{i=1}^d [S^{k_i,p_j}_j]}{([S^{p_j}_j])^{d-1}} 
	~ \Big| ~ d_j = d     \bigg] ,\qquad d=2,3,\dots
	\label{eq:S_B_d-funnel-equality}
\end{equation}
\end{lemma}
\begin{proof}
% The case $d=1$ is immediate. 
%	
%If~$j$ does not lie on any cycle, then 
%by the funnel equality
%	   \begin{subequations}
%	   	\label{eqs:funnel-ineq-ER-d}
%
%
%	If, however, $j$ lies on $N_j$ cycles $\{C_n\}_{n=1}^{N_j}$, such that the number of nodes 
%	of cycle~$C_n$ is~$L_n$, then  by Theorem~\ref{thm:upper-bound-funnel-d},
%	%by~\eqref{eqs:funnel-upper-bound-Ncycle}, 
%	and denoting $\tilde{q}:=\frac{q}{\lambda}$, we have
%	\begin{subequations}
	%		\label{eqs:funnel-ineq-ER-d-upper-bound}
	%		\begin{equation}
		%		0<	[S_j] -
		%			\frac{\prod_{i=1}^d [S^{k_i,p_j}_j]}{([S^{p_j}_j])^{d-1}} 
		%			\le  2 [S_j^{p_j}] 
		%			\sum_{n=1}^{N_j}E(t;L_n) ,
		%		\end{equation}
	%		where
	%		\begin{equation}
		%			E(t;L_n) \le e^{-\left(p+\tilde{q}\right)t} \left(\frac{e\tilde{q}t}{\left \lfloor{\frac{L_n+1}{2}}\right \rfloor}\right)^{\left \lfloor{\frac{L_n+1}{2}}\right \rfloor}, 
		%			\qquad    0<t< \frac1{\tilde q}  \bigg \lfloor{\frac{L_n+1}{2}} \bigg \rfloor, 
		%			%\qquad  \left \lfloor{\frac{L_n+1}{2}}\right \rfloor \ge \tilde{q}t, 
		%			%	\qquad 
		%			%	E(t;L_n) \le  \left(\frac{q}{p+q}\right \rfloor}\right)^{\left \lfloor{\frac{L_n+1}{2}}\right \rfloor},
	%		%	\qquad   t \ge 0.
	%	\end{equation}
%	and
%	\begin{equation}
	%		%	\label{eq:funnel-ineq-ER-d-upper-bound}
	%		E(t;L_n) \le  \bigg(\frac{\tilde{q}}{p+\tilde{q}} \bigg)^{\left \lfloor{\frac{L_n+1}{2}}\right \rfloor},
	%		\qquad   0<t .
	%	\end{equation}
%\end{subequations}
Let $j \in B_{d}$ and denote $\tilde{q}:=\frac{q}{\lambda}$.  
Fix $t>0$, and let $L^0 = L^0(t):= \left \lceil 2\tilde{q}t-1 \right \rceil$.
By Theorem~\ref{thm:funnel-ER}, 
and since $[S^{p_j}_j] \le 1$, see~\eqref{eq:[S^pj_j](t;N^ER)=(1-beta0)e^-pt}, 
$$
\begin{aligned}
	0 \le 
	[S_j] -	\frac{\prod_{i=1}^d [S^{k_i,p_j}_j]}{([S^{p_j}_j])^{d-1}}  
	%\\  &\
	\le
	2 
	\sum_{L=3}^{L^0(t)}c_L(j)
	\bigg(\frac{\tilde{q}}{p+\tilde{q}} \bigg)^{\left \lfloor{\frac{L+1}{2}}\right \rfloor} 
	+ 	2\sum_{L=L^0(t)+1}^{M}c_L(j) e^{-\left(p+\tilde{q}\right)t}
	\bigg(\frac{e\tilde{q}t}{\left \lfloor{\frac{L+1}{2}}\right \rfloor}\bigg)^{\left \lfloor{\frac{L+1}{2}}\right \rfloor},
\end{aligned}
$$
where $c_L(j)$ is the number of cycles of length~$L$ that go through~$j$.
Taking the expectation of these inequalities,  
conditional on~$d_j=d$, 
%Multiplying these inequalities by $\mathbbm{1}_{d_j=d}$
%	Averaging these inequalities over all $j \in B_{d}$	
%averaging over~${\cal G}(M,\frac{\lambda}{M})$,
%	noting that
%	$$
%	\mathbb{E}_{{\cal G}(M,\frac{\lambda}{M})} \Big[ c_L(j) ~ \big| ~ d_j = d  \Big]
%	=
%	\mathbb{E}_{{\cal G}(M,\frac{\lambda}{M})} \Big[\frac{1}{|B_{d}|}\sum_{m\in B_{d}}c_L(m;G)\Big]
%	=
%	\mathbb{E}_{{\cal G}(M,\frac{\lambda}{M})}	\big[c_L^d\big],
%	$$
and using~\eqref{eq:c_L^d}, gives
\begin{subequations}
	\label{eqs:E_G[S_B_d]}
	\begin{equation}
		\begin{aligned}
			0 &
			\le 
			\mathbb{E}_{{\cal G}(M,\frac{\lambda}{M})} \bigg[ [S_j] -	  \frac{\prod_{i=1}^d [S^{k_i,p_j}_j]}{([S^{p_j}_j])^{d-1}} 
			~ \Big| ~ d_j = d     \bigg]
			\\  &\le
			2  \sum_{L=3}^{L_0(t)} C_L^d \, 
			\bigg(\frac{\tilde{q}}{p+\tilde{q}} \bigg)^{\left \lfloor{\frac{L+2}{2}}\right \rfloor} 
			+ 
			2\sum_{L=L_0(t)+1}^{M} C_L^d \,  %e^{-\left(p+q\right)t}
			\bigg(\frac{e\tilde{q}t}{\left \lfloor{\frac{L+2}{2}}\right \rfloor}\bigg)^{\left \lfloor{\frac{L+2}{2}}\right \rfloor}.
		\end{aligned}
	\end{equation}
	%where $c_L^d$ is given by~\eqref{eq:c_L^d}. 
	%	where  $c_L(d,G^{\rm ER})$ is the number of different cycles with $L$~nodes that go through a node of degree~$d$.
	%	Recall that, see~\eqref{eq:E_G[c_L^d]},
	%		\begin{equation*}
		%			C_L^d \sim \frac{ d (d-1)}{2} \lambda^{L-2} \frac1{M}, \qquad M \to \infty.
		%	\end{equation*}
	Using the upper bound~\eqref{eq:E_G[c_L^d]-bound} for~$C_L^d$, we have
	\begin{equation}
		\sum_{L=3}^{L_0(t)} C_L^d \, 
		\bigg(\frac{\tilde{q}}{p+\tilde{q}} \bigg)^{\left \lfloor{\frac{L+2}{2}}\right \rfloor}
		<
		\sum_{L=3}^{L_0(t)}\frac{\kappa_{L,d}^{\rm ER}}{M}
		\bigg(\frac{\tilde{q}}{p+\tilde{q}} \bigg)^{\left \lfloor{\frac{L+2}{2}}\right \rfloor}
		= O\left(\frac1M \right), \quad M \to \infty,
	\end{equation}
	where $\kappa_{L,d}^{\rm ER}:=\frac{ d (d-1)}{2} \lambda^{L-2}$, and 
	\begin{equation}
		\begin{aligned}
			&	\sum_{L=L_0(t)+1}^{M} C_L^d \, 
			\bigg(\frac{e\tilde{q}t}{\left \lfloor{\frac{L+2}{2}}\right \rfloor}\bigg)^{\left \lfloor{\frac{L+2}{2}}\right \rfloor}
			<
			\sum_{L=L_0(t)+1}^{M} \frac{\kappa_{L,d}^{\rm ER}}{M}
			\bigg(\frac{e\tilde{q}t}{\left \lfloor{\frac{L+2}{2}}\right \rfloor}\bigg)^{\left \lfloor{\frac{L+2}{2}}\right \rfloor}
			\\		\vspace{20mm}
			& \quad < 
			\frac1{M}  \frac{ d (d-1)}{2} \sum_{L=L_0(t)+1}^{\infty}  \lambda^{L-2} 
			\bigg(\frac{e\tilde{q}t}{\left \lfloor{\frac{L+2}{2}}\right \rfloor}\bigg)^{\left \lfloor{\frac{L+2}{2}}\right \rfloor}
			= O\left(\frac1M \right), \qquad M \to \infty.
		\end{aligned}
	\end{equation}
\end{subequations}
%Since $	\overline{\left[S^\infty_{B_{d}}\right]}:= 
%\lim_{M \to \infty} \mathbb{E}_{{\cal G}(M,\frac{\lambda}{M})} \left[\, [S_{B_{d}}] \, \right]$, see, 
Letting $M \to \infty$   in relations~\eqref{eqs:E_G[S_B_d]}
and using~\eqref{eq:E-E-R-j} and~\eqref{eq:E-E-R-M=infty}  proves~\eqref{eq:S_B_d-funnel-equality}.
\end{proof}

When $\lambda>1$, the number of cycles of length~$L$ increases exponentially as~$\lambda^L$,
see~\eqref{eq:E_G[c_L^d]}.
By~\eqref{eq:funnel-upper-bound-Ncycle},
%each node belongs, on average, to at least 
%$3 \left(\frac{\lambda}{2}-1 \right)$ cycles (see Lemma~\ref{lem:c^ER>3(lambda/2-1)}), 
a strict funnel inequality holds for any node that lies on a cycle.
Nevertheless,
Lemma~\ref{lem:S_B_d-funnel} shows that one can 
compute~$\overline{\left[S^\infty_{B_{d}}\right]}$ using the funnel {\em equality}~\eqref{eq:funnel-ineq-ER-d}. 
This is because as $L$ increases, 
%the number of cycles increases exponentially as~$\lambda^L$, see~\eqref{eq:E_G[c_L^d]}, but  
the effect of cycles on the adoption probability decays 
at the super-exponential rate of~$\frac{1}{L^L}$,
see~\eqref{eq:funnel-upper-bound-Ncycle}. Therefore, 
{\em the overall effect of the cycles vanishes as $M \to \infty$}.

\subsection{$\protect\overline{\left[S^\infty_{B_{d}}\right]}=\left(\frac1{e^{-pt}(1-I^0)}\right)^{d-1}
	\left(\overline{\big[S^\infty_{B_{1}}\big]}\right)^{d}$}

Next, we express the expected susceptibility level of nodes of degree~$d$ with 
that of nodes of degree~$1$: 
\begin{lemma}
	\label{lem:lemma3}
	%	Consider the Bass model~\eqref{eqs:Bass-ER} on $ \lim_{M \to \infty} {\cal G}(M,\frac{\lambda}{M})$. Then
	Let $\overline{\left[S^\infty_{B_{d}}\right]}$ be given by~{\em (\ref{eq:S_Bd_by_dist},\ref{eq:E-E-R-M=infty})}. Then 
	\begin{equation}
		\overline{\left[S^\infty_{B_{d}}\right]}=\left(\frac1{e^{-pt}(1-I^0)}\right)^{d-1}
		\left(\overline{\big[S^\infty_{B_{1}}\big]}\right)^{d},\qquad d=1,2,\dots
		\label{eq:s_k-ER}
	\end{equation}
\end{lemma}
\begin{proof}
	The case $d=1$ is immediate.  For $d \ge 2$,
	we compute the right-hand side of~\eqref{eq:S_B_d-funnel-equality}
	as follows. 
	Since $[{S_{{j}}^{p_{{j}}}}]= (1-I^0)e^{-pt}$ is independent of~${{\cal N}}^{{\rm ER}}(G)$,
	see~\eqref{eq:[S^pj_j](t;N^ER)=(1-beta0)e^-pt}, 
	\begin{subequations}
		\label{eqs:end-pf-ER}
		\begin{equation}
			\mathbb{E}_{{\cal G}} \bigg[\, 
			\frac{\prod_{i=1}^d [S^{k_i,p_j}_j]}{([S^{p_j}_j])^{d-1}}  ~ \Big| ~ d_j = d     \bigg]
			= 
			\left(\frac{1}{e^{-pt} (1-I^0)}\right)^{d-1}	
			\mathbb{E}_{{\cal G}} 
			\bigg[	 \prod_{i=1}^d [S^{k_i,p_j}_j]  ~ \Big| ~ d_j = d    \bigg].
		\end{equation}
		%In addition, $\mathbb{E}_{{\cal G}(M,\frac{\lambda}{M})} \left[\prod_{i=1}^d  [S^{k_i,p_j}_j] \right] $ is independent of~$j$.
		%Therefore, for any ${j} \in \cal M$, 
		%	\begin{equation}
			%\frac{1}{M}\sum_{j =1}^M \mathbb{E}_{{\cal G}(M,\frac{\lambda}{M})} \bigg[ \prod_{i=1}^d 	 [S^{k_i,p_j}_j] \bigg]
			%= 
			%	\mathbb{E}_{{\cal G}(M,\frac{\lambda}{M})} \bigg[ \prod_{i=1}^d  [S^{k_i,p_{{j}}}_{{j}}] \bigg] .
			%\end{equation}
			By the indifference principle, we can compute 
			$\{[S_{{j}}^{A_{i}, p_{{j}}}](t;{\cal N}^{{\rm ER}})\}$ on an equivalent network, denoted by~$\widetilde{{\cal N}}^{{\rm ER}}$,
			in which we delete the $d$~directional edges $\{{j} \to k_i\}_{i=1}^d$, but leave the  $d$~directional edges $\{k_i \to {j}\}_{i=1}^d$ intact. Thus, 
			\begin{equation}
				[S_{{j}}^{A_{i}, p_{{j}}}] = [\widetilde{S_{{j}}^{A_{i}, p_{{j}}}}],
				\qquad i=1, \dots, d,
			\end{equation}
			where $[\widetilde{\, \cdot \,} ]$ denote nonadoption probabilities in~$\widetilde{{\cal N}^{{\rm ER}}}$. 
			By~\cite[Lemma~7.2]{Funnel-25},  
			$[\widetilde{S_{{j}}^{A_{i}, p_{{j}}}}]  =
			[\widetilde{S_{{j}}^{A_{i}}}] \,
			[\widetilde{S_{{j}}^{p_{{j}}}}]$,
			where $[\widetilde{S_{{j}}^{A_{i}}}]:=[\widetilde{S_{{j}}^{A_{i}, p_{{j}}:=0}}]$.
			Since $[\widetilde{S_{{j}}^{p_{{j}}}}]= (1-I^0)e^{-pt}$ is independent of~$\widetilde{{\cal N}}^{{\rm ER}}$,
			\begin{equation}
				\begin{aligned}
					\mathbb{E}_{{\cal G}} 
					\bigg[	\prod_{i=1}^d [S^{k_i,p_j}_j]   ~ \Big| ~ d_j = d   \bigg]
					&=	\mathbb{E}_{{\cal G}} 
					\bigg[	 \prod_{i=1}^d [\widetilde{S_{{j}}^{A_{i}, p_{{j}}}}]   ~ \Big| ~ d_j = d   \bigg] 			
					%\\ &
						= 
					\Big((1-I^0)e^{-pt}\Big)^d \,
					\mathbb{E}_{{\cal G}} 	 			
					\bigg[  \prod_{i=1}^d  [\widetilde{S_{{j}}^{A_{i}}}]
					~ \Big| ~ d_j = d  \bigg]. 
				\end{aligned}
			\end{equation}
			
			The nonadoption probability~$[\widetilde{S_{{j}}^{A_{i}}}]$ is evaluated 
			in the network~$\widetilde{{\cal N}}^{{\rm ER},A_i}
			:=\widetilde{{\cal N}}^{{\rm ER},A_i,p_{{j}}:=0}$.
			In this network,  
			%$\widetilde{{\cal N}}^{A_i}:=\widetilde{{\cal N}}^{A_i,p_{{j}}:=0}$, 
			${j}$ is a degree-one node that can only adopt because of internal influences by the node~$k_i$.  
			%by the indifference principle, we can delete the edge ${j} \to k_i$ when
			%computing $[S^{k_i,p_{{j}}}_{{j}}]$. Therefore, 
			%$[S^{k_i,p_{{j}}}_{{j}}]$ is completely determined by~$[S_{k_i}]$,
			%	see~\eqref{eq:d_dt_[Sj]-i--j-general}. 
			Since the only thing that relates the  nodes $\{k_{i}\}_{i=1}^d$ 
			in the network~$\widetilde{{\cal N}}^{{\rm ER}}$ is that they are connected to~${j}$, the conditional  probabilities~$\big\{\big[\widetilde{S_{k_{i}}}\big] \, \big| \, d_j = d  \big\}_{i=1}^d$
			become independent in~${\cal G}(M,\frac{\lambda}{M})$ as $M \to \infty$.
			Therefore, 
			%the nonadoption probabilities 
			$\big\{ \big[\widetilde{S_{{j}}^{A_{i}}} \big] \, \big| \, d_j = d  \big\}_{i=1}^d$
			become independent in~$\lim_{M \to \infty} {\cal G}(M,\frac{\lambda}{M})$ as well.  
			Hence,
			\begin{equation}
				\label{eq:end-pf-d-regular-product}
				\lim_{M \to \infty} 
				%	 \frac{1}{M}\sum_{j =1}^M 
				\mathbb{E}_{{\cal G}(M,\frac{\lambda}{M})} \bigg[    \prod_{i=1}^d  \big[\widetilde{S_{{j}}^{A_{i}}}\big] 
				~ \Big| ~ d_j = d  \bigg]
				=
				%	 \frac{1}{M}\sum_{j =1}^M 
				\prod_{i=1}^d 
				\lim_{M \to \infty} 
				\mathbb{E}_{{\cal G}(M,\frac{\lambda}{M})} \left[    \big[\widetilde{S_{{j}}^{A_{i}}}\big]  ~ \Big| ~ d_j = d   \right].
			\end{equation}
			% 		{\bf Therefore},
			%	 		\begin{equation}
				%	\label{eq:end-pf-d-regular-product}
				%	\lim_{M \to \infty} 
				%	%	 \frac{1}{M}\sum_{j =1}^M 
				%	\mathbb{E}_{{\cal G}(M,\frac{\lambda}{M})} \bigg[	 \mathbbm{1}_{d_j=d} \,  \prod_{i=1}^d  [\widetilde{S_{{j}}^{A_{i}}}] \bigg]
				%	=
				%	%	 \frac{1}{M}\sum_{j =1}^M 
				%	\prod_{i=1}^d 
				%	\lim_{M \to \infty} 
				%	\mathbb{E}_{{\cal G}(M,\frac{\lambda}{M})} \bigg[  	 \mathbbm{1}_{d_j=d} \,  [\widetilde{S_{{j}}^{A_{i}}}]  \bigg].
				%\end{equation}
				Using again the relations $[\widetilde{S_{{j}}^{A_{i}, p_{{j}}}}]  =
				[\widetilde{S_{{j}}^{A_{i}}}] \,
				[\widetilde{S_{{j}}^{p_{{j}}}}] = [\widetilde{S_{{j}}^{A_{i}}}] \,(1-I^0)e^{-pt}$,  we have
				\begin{equation}
					\begin{aligned}
						&	\Big((1-I^0)e^{-pt}\Big)^d \,	\prod_{i=1}^d 
						\lim_{M \to \infty} 
						\mathbb{E}_{{\cal G}(M,\frac{\lambda}{M})} 
						\bigg[ \big[\widetilde{S^{k_i}_{{j}}}\big]  ~ \Big| ~ d_j = d   \bigg]
						=
						\prod_{i=1}^d 
						\lim_{M \to \infty} 
						\mathbb{E}_{{\cal G}(M,\frac{\lambda}{M})} \bigg[  \big[\widetilde{S^{k_i,p_{{j}}}_{{j}}}\big]
						~ \Big| ~ d_j = d    \bigg]
						\\ \quad & =
						\prod_{i=1}^d 
						\lim_{M \to \infty} 
						\mathbb{E}_{{\cal G}(M,\frac{\lambda}{M})} \bigg[  \big[{S^{k_i,p_{{j}}}_{{j}}}\big] 
						~ \Big| ~ d_j = d   \bigg].
					\end{aligned}
				\end{equation}
				%	 	\end{subequations}
			%	 	The result follows from relations~\eqref{eqs:end-pf-d-reg}. 
			%	 
			%	 
			%	  By~\eqref{eq:[S^pj_j](t;N^ER)=(1-beta0)e^-pt}, 
			% \begin{subequations}
				% 	\label{eqs:end-pf-ER}
				%	\begin{equation}
					%	\mathbb{E}_{G} \bigg[\, 
					%	 \mathbbm{1}_{d_j=d} \, \frac{\prod_{i=1}^d [S^{k_i,p_j}_j]}{([S^{p_j}_j])^{d-1}}    \bigg] 
					%	= 
					%	\left(\frac1{e^{-pt}(1-I^0)}\right)^{d-1}
					%	\mathbb{E}_{G} \bigg[\, 
					%	 \mathbbm{1}_{d_j=d} \, \prod_{i=1}^d [S^{k_i,p_j}_j]    \bigg] .
					%\end{equation}
					%
					%	  Since $[S^{k_i,p_j}_j]$ is completely determined by~$[S_{k_i}]$,
					%see~\eqref{eq:d_dt_[Sj]-i--j-general},
					%and since~$\{[S_{k_{1}}], \dots, [S_{k_{d}}] \}$
					%become independent as $M \to \infty$, 
					%the probabilities $\left\{[S^{k_i,p_j}_j](t;{\cal N}^{\rm ER})\right\}_{i=1}^d$
					%become independent in~$\lim_{M \to \infty} {\cal G}(M,\frac{\lambda}{M})$ as well.  
					%Hence,
					%	\begin{equation}
						%	\lim_{M \to \infty} 
						%	\mathbb{E}_{\cal G} \bigg[\, 
						%	 \mathbbm{1}_{d_j=d} \, \prod_{i=1}^d [S^{k_i,p_j}_j]    \bigg] 
						%	=
						%	\prod_{i=1}^d	\lim_{M \to \infty} 
						%	\mathbb{E}_{\cal G} \bigg[\, 
						%	 \mathbbm{1}_{d_j=d} \, [S^{k_i,p_j}_j]    \bigg] .
						%\end{equation}
						Finally, since~${\cal N}^{A_i,p_j}$  is an ER network in ${\cal G}(M,\frac{\lambda}{M})$ in which $j$~is a degree-one node,
						%	 and it can only be influenced by nodes in~$A_i$. Therefore,
						%	 $j$~is a degree-one node in an  ER network in $G\left(\frac{M}{d},\alpha\right)$
						\begin{equation}
							\label{eq:E_G[[S^A_i,p_j_j]=[S_B_1]}
							\lim_{M \to \infty} 
							\mathbb{E}_{\cal G}\bigg[\, 
							[S^{k_i,p_j}_j]  ~ \Big| ~ d_j = d  \bigg] 
							= \overline{\big[S^\infty_{B_{1}}\big]}, \qquad i=1, \dots, d.
						\end{equation} 
					\end{subequations}
					%	  Hence, using~\eqref{eq:E_G[[S^A_i,p_j_j]=[S_B_1]}, 
					%	  \begin{equation}
						%	  	\label{eq:lim_M->infty_E_G[product]}
						%	  \lim_{M\to \infty}	\mathbb{E}_G\left[\prod_{i=1}^d [S^{k_i,p_j}_j] \right] = 
						%	  	\prod_{i=1}^d \lim_{M\to \infty} \mathbb{E}_G\left[ [S^{k_i,p_j}_j] \right]
						%	  	 =  \left(\overline{\big[S^\infty_{B_{1}}\big]}\right)^{d}.
						%	  \end{equation}
					The result follows from~\eqref{eq:S_B_d-funnel-equality}
					and~\eqref{eqs:end-pf-ER}.
					\end{proof}

	\subsection{Degree-one nodes}
	
	So far, we saw that, see~\eqref{eq:s_by_dist} and~\eqref{eq:s_k-ER}, 
	$$
	\overline{[S^\infty]}=\sum_{d=0}^{\infty}\mathbb{P}_{\rm deg}(d) \,
	%	\lim_{M\rightarrow\infty}
	\overline{[S^\infty_{B_{d}}]}, \qquad 
	\overline{\left[S^\infty_{B_{d}}\right]}=\left(\frac1{e^{-pt}(1-I^0)}\right)^{d-1}
	\left(\overline{\big[S^\infty_{B_{1}}\big]}\right)^{d}.
	$$
	Therefore, the 
	calculation of $\overline{[S^\infty]}$ reduces to that of~$\big[S^\infty_{B_{1}}\big]$.  
	To compute~$\big[S^\infty_{B_{1}}\big]$, 
	we distinguish between nodes of degree one according to the degree of the node to which they are connected.
	Thus, let 
	
	\begin{equation}
		\label{eq:S_B1_d:ER}
		%	\overline{\left[S^\infty_{B_{d}}\right]}:= \lim_{M \to \infty}\overline{[S_{B_{d}}]}, \qquad 
		\overline{[S_{B_{1,d}}]}:=\mathbb{E}_{{\cal G}(M,\frac{\lambda}{M})}\Big [ \,  [S_{j}]~ \Big| ~ d_j = 1,
		d_{k_{1}(j)}=d   \Big],
	\end{equation}
	where $k_{1}(j)$ is the unique node connected to~$j$. Thus, $\overline{[S_{B_{1,d}}]}$ is the expected susceptibility level of nodes of degree one that are connected to a node of degree~$d$. 
	%
	%
	%
	%
	%$$
	%B_{1,d}(G):=\big\{ m\in B_1  ~ \big| ~ e_{m,m_{1}}=1,d_{m_{1}}=d \big\}
	%$$ 
	%%Denote 
	%denote the set of all the nodes of degree~$1$ that are connected to nodes of degree~$d$
	%in~$G$.
	%Denote also 
	%%\[
	%%B_{1,d}(G):=\{ m\in{\cal M} \mid d_{m}=1,e_{m,m_{1}}=1,d_{m_{1}}=d \}, \qquad d=1,2,\dots,
	%%\]
	%the expected nonadoption level in~$B_{1,d}$ by
	%\begin{equation}
	%	[S_{B_{1,d}}]\big(t;{\cal N}^{\rm ER}(G)\big)
	%	:=\frac{1}{\left|B_{1,d}\right|}\sum_{m\in B_{1,d}}[S_m],
	%	\label{eq:sa1i}
	%\end{equation}
	%and the average of $[S_{B_{1,d}}]$ over ${\cal G}^{\rm ER}(M,\frac{\lambda}{M})$ by
	%\begin{equation}
	%	\overline{[S_{B_{1,d}}]}(t):=\mathbb{E}_{\cal G}[S_{B_{1,d}}].
	%	\label{eq:sali2}
	%\end{equation}
	
	\begin{lemma}
		\label{lem:lemma2-1}
		%			Consider the Bass model~\eqref{eqs:Bass-ER} on $\lim_{M \to \infty} {\cal G}(M,\frac{\lambda}{M})$.			Then
		Let
		\begin{equation}
			\overline{\big[S^\infty_{B_{1,d}}\big]}(t):= \lim_{M \to \infty} \overline{[S_{B_{1,d}}]}(t),
			\qquad d=1,2, \dots
			\label{eq:sali2-M=infty}
		\end{equation}
		Then
		\begin{equation}
			\overline{\big[S^\infty_{B_{1}}\big]}=\sum_{d=1}^{\infty}\mathbb{P}_{\rm deg}(d-1)\overline{\big[S^\infty_{B_{1,d}}\big]},\label{eq:averaging}
		\end{equation}
		where $\mathbb{P}_{\rm deg}$ is  given by \eqref{eq:pois}.
		% and $\overline{\big[S^\infty_{B_{1,d}}\big]}$
		%	are given by \eqref{eq:pois} and \eqref{eq:sali2-M=infty}, respectively.
	\end{lemma}
	\begin{proof}
		By the law of total expectation,  \eqref{eq:S_Bd_by_dist}, and~\eqref{eq:S_B1_d:ER},
		$$
		\begin{aligned}
			\overline{[S_{B_1}]}
			& = \mathbb{E}_{{\cal G}(M,\frac{\lambda}{M})}\big [ \,  [S_{j}]~ \big| ~ d_j = 1  \big]
			=\sum_{d=0}^{M-1} 
			\mathbb{P}\big(d_{k_1}(j) = d ~ \big| ~d_j = 1 \big) \,    %\mathbb{P}_{\rm deg}^{M}(d)
			\mathbb{E}_{{\cal G}(M,\frac{\lambda}{M})}\Big [ \,  [S_{j}]~ \Big| ~ d_j = 1,
			d_{k_{1}(j)}=d   \Big]
			\\ &= \sum_{d=0}^{M-1} 
			\mathbb{P}\big(d_{k_1}(j) = d ~ \big| ~d_j = 1 \big) \, \overline{[S_{B_{1,d}}]}.
		\end{aligned}
		$$	 
		Since  
		$\mathbb{P}\big(d_{k_1}(j) = d ~ \big| ~d_j = 1 \big)$ is  the probability that
		the node~$k_{1}$ has $d-1$ other edges in a network with $M-1$ ~nodes, 
		\begin{equation*}
			\mathbb{P}\big(d_{k_1}(j) = d ~ \big| ~d_j = 1 \big) = \mathbb{P}_{\rm deg}^{M-1}(d-1),
			\qquad d=1,\dots, M-1.
			%\label{eq:spam}
		\end{equation*}
		%	Hence,
		%	$$
		%\overline{\big[Sy_{B_{1}}\big]}
		%=\sum_{d=0}^{M-1} 
		%\mathbb{P}_{\rm deg}^{M-1}(d-1) \,    %\mathbb{P}_{\rm deg}^{M}(d)
		%\overline{\big[S_{B_{1,d}}\big]}.
		%$$	 
		Letting $M \to \infty$ and using~\eqref{eq:pois},
		\eqref{eq:E-E-R-M=infty}, and~\eqref{eq:sali2-M=infty} yields the result. 
	\end{proof}

\subsection{Master equation for $\protect\overline{\big[S^\infty_{B_{1,d}}\big]}$}

Next, we derive the master equation for~$\big[S^\infty_{B_{1,d}}\big]$:
\begin{lemma}
	\label{lem:master_s_1_k}
	Let $\big[S^\infty_{B_{1,d}}\big]$ be given by~\eqref{eq:sali2-M=infty}.
	%where $d  \in \mathbb{N}$. 
	Then 
	%the master equation for~$\big[S^\infty_{B_{1,d}}\big]$ is given by
	%\begin{equation}
	%	\overline{\big[S^\infty_{B_{1,d}}\big]}(t):= \lim_{M \to \infty} \overline{[S_{B_{1,d}}]}(t),
	%	\qquad d=1,2, \dots
	%	\label{eq:sali2-M=infty}
	%\end{equation}
	%		
	%		
	%		Consider the Bass model~\eqref{eqs:Bass-ER} on $\lim_{M \to \infty} {\cal G}(M,\frac{\lambda}{M})$.
	%\begin{subequations}
	\begin{equation}
		\frac{d \overline{\big[S^\infty_{B_{1,d}}\big]}}{dt}=-\left(p+\frac{q}{\lambda}\right)\overline{\big[S^\infty_{B_{1,d}}\big]}
		+\frac{q}{\lambda}(1-I^0)e^{-pt}\overline{\left[S^\infty_{B_{d-1}}\right]},
		\quad t>0, \qquad 	\overline{\big[S^\infty_{B_{1,d}}\big]}(0)=1-I^0,
		\label{eq:master_s_1_k}
	\end{equation}
	where $\overline{\left[S^\infty_{B_{d-1}}\right]}$ is given by~\eqref{eq:E-E-R-M=infty}.
	%  subject to the initial condition
	%		\begin{equation}
		%			\overline{\big[S^\infty_{B_{1,d}}\big]}(0)=1-I^0.
		%			\label{eq:init7}
		%		\end{equation}
	%\end{subequations}
\end{lemma}
 \begin{proof}
	The initial condition follows from~\eqref{eq:Bass-ER-IC}.
	To derive the ODE, recall that the master equation for~$[S_j]$ on a general network 
	is given by, see~\cite[Lemma~2]{MOR-22}, 
	\[
	\frac{d[S_j]}{dt}=-\Big(p_m+\sum_{k \in \cal M}q_{k,j}\Big)[S_j]+\sum_{k \in \cal M}q_{k,j}[S_{k,j}], \qquad 
	[S_j](0)=1-I^0.
	\]
	%	Let ${\cal N}^{\rm ER} = {\cal N}^{\rm ER}(G)$, where $G \in {\cal G}^{\rm ER}(M,\frac{\lambda}{M})$.
	%d_{k_{1}(j)}=d $.
	Substituting $d_j = 1$ and the ER~network structure~\eqref{eq:ER-network}  gives
	\begin{equation}
		\frac{d [S_j]}{dt}=-\left(p+\frac{q}{\lambda}\right)[S_j]+\frac{q}{\lambda}[S_{k_1(j),j}],
		\label{eq:master2}
	\end{equation}
where $k_{1}(j)$ is the unique node connected to~$j$. By the {\em indifference principle}
\cite[Lemma~3.8]{Bass-boundary-18}, 
	%such that $\text{\rm degree}(k_{1}) = d$.
%	By the indifference principle (Theorem~\ref{thm:Indifference}), %\cite[lemma~3.8]{2}, 
	$[S_{k_1,j}]$ can be calculated on an equivalent network~$\widetilde{\cal N}$,
	in which the non-influential edge  $j \leftrightarrow k_{1}$
	is removed, i.e.,
	\begin{equation}
		[S_{k_1,j}]=\widetilde{[S_{k_1,j}],}
		\label{eq:smsm1}
	\end{equation}
	where the tildes denote probabilities in $\widetilde{\cal N}$.
	Since  $d_j=1$ in the original network, $j$~is an isolated node
	in~$\widetilde{\cal N}$. Hence, $\widetilde{[S_j]}=(1-I^0)e^{-pt}$, 
	see~\eqref{eq:[S^pj_j](t;N^d-regular)=(1-beta0)e^-pt}, and so 
	\begin{equation}
		\widetilde{[S_{k_1,j}]}=\widetilde{[S_j]}\widetilde{\left[S_{k_{1}}\right]}=(1-I^0)e^{-pt}\widetilde{\left[S_{k_{1}}\right]}.
		\label{eq:indifference}
	\end{equation}
	Substituting (\ref{eq:smsm1}) and~(\ref{eq:indifference})
	in~(\ref{eq:master2}), we get
	\begin{equation}
		\frac{d [S_j]}{dt}=-\left(p+\frac{q}{\lambda}\right)[S_j]+\frac{q}{\lambda}(1-I^0)e^{-pt}
		\widetilde{\left[S_{k_{1}}\right]}.\label{eq:ddtsm}
	\end{equation}
	Taking the expectation of~\eqref{eq:ddtsm} over~${\cal G}(M,\frac{\lambda}{M})$,
	conditional on $d_j=1$ and $d_{k_1(j)}=d$, 
	gives
	\begin{equation}
		\frac{d [S_{B_{1,d}}]}{dt}=
		-\left(p+\frac{q}{\lambda}\right)[S_{B_{1,d}}]+
		\frac{q}{\lambda}(1-I^0)e^{-pt}	\,
		\mathbb{E}_{ {\cal G}(M,\frac{\lambda}{M})} \left[\widetilde{\left[S_{k_{1}(j)}\right]} ~ \Big| ~ d_j = 1, d_{k_1(j)}=d \right]
		.\label{eq:egsak1}
	\end{equation}
	The nonadoption probability $\Big[\widetilde{\left[S_{k_{1}(j)}\right]} ~ \Big| ~ d_j = 1, d_{k_1(j)}=d \Big]$
	can be considered on the network~$\widetilde{\widetilde{\cal N}}$,
	which is~$\widetilde{\cal N}$ without the isolated node~$j$. Network $\widetilde{\widetilde{\cal N}}$ 
	is thus in~$ G\left(M-1,\frac{\lambda}{M}\right)$, such that $d_{k_{1}} = d-1$ in~$\widetilde{\widetilde{\cal N}}$. Hence,
	\begin{equation}
		%	\mathbb{E}_{ {\cal G}(M,\frac{\lambda}{M})}\left[\frac{1}{\left|B_{1,d}\right|}\sum_{m\in %B_{1,d}}\widetilde{\left[S_{k_{1}(m)}\right]}\right](t;M)
		%
		\mathbb{E}_{ {\cal G}(M,\frac{\lambda}{M})} \left[\widetilde{\left[S_{k_{1}(j)}\right]} ~ \Big| ~ d_j = 1, d_{k_1(j)}=d \right]
		=\overline{\left[S_{B_{d-1}}\right]}(t;M-1).\label{eq:egsak2}
	\end{equation}
	Combining (\ref{eq:egsak1}) and (\ref{eq:egsak2}) and
	letting $M\rightarrow\infty$ gives (\ref{eq:master_s_1_k}).
 \end{proof}

\subsection{Solving the master equations}

We are now in a position to solve the master equations.

\begin{lemma}
	\label{lem:lemma4}
%	Consider theBass and {\rm SI} models~\eqref{eqs:Bass-ER} on ${\cal G}^{\rm ER}(M,\frac{\lambda}{M})$.
	Let 
	\begin{equation}
		\label{eq:y=e^(pt)(1-beta_0)[S^infty_B_1]}
		y(t):=\frac{\overline{\big[S^\infty_{B_{1}}\big]}}{e^{-pt}(1-I^0)}.
	\end{equation}
	Then
	\begin{enumerate}
		\item $y(t)$ is the solution of equation~\eqref{eq:y-ER}. 
		\item 
		$
		f^{\rm ER}(t) = 1-(1-I^0) e^{-\lambda  -pt +\lambda y(t)}
		$. 
	\end{enumerate}
\end{lemma}
 \begin{proof}
 	Since $\big[S^\infty_{B_{1}}\big](0)= 1-I^0$, see~\eqref{eq:Bass-ER-IC}, then
 $y(0) = 1$.  
 By~\eqref{eq:averaging}, 
 \[
 \frac{d \overline{\big[S^\infty_{B_{1}}\big]}}{dt}=\sum_{d=1}^{\infty}\mathbb{P}_{\rm deg}(d-1)\frac{d \overline{\big[S^\infty_{B_{1,d}}\big]}}{dt}.
 \]
 In addition, substituting (\ref{eq:s_k-ER}) in~(\ref{eq:master_s_1_k}) gives
 \begin{equation*}
 	\frac{d \overline{\big[S^\infty_{B_{1,d}}\big]}}{dt}=-\left(p+\frac{q}{\lambda}\right)\overline{\big[S^\infty_{B_{1,d}}\big]}+\frac{q}{\lambda}\left(\frac1{e^{-pt}(1-I^0)}\right)^{d-3}\left(\overline{\big[S^\infty_{B_{1}}\big]}\right)^{d-1},\qquad d=1,2,\dots
 	%\label{eq:funnel}
 \end{equation*}
 Combining these two relations gives 
 \[
 \frac{d \overline{\big[S^\infty_{B_{1}}\big]}}{dt}
 =\sum_{d=1}^{\infty}\mathbb{P}_{\rm deg}(d-1)\left[-\left(p+\frac{q}{\lambda}\right)\overline{\big[S^\infty_{B_{1,d}}\big]}+\frac{q}{\lambda}\left(\frac1{e^{-pt}(1-I^0)}\right)^{d-3}\left(\overline{\big[S^\infty_{B_{1}}\big]}\right)^{d-1}\right].
 \]
 Therefore, using~\eqref{eq:pois} and~(\ref{eq:averaging}), we get
 \begin{equation}
 	\begin{aligned}
 		\frac{d \overline{\big[S^\infty_{B_{1}}\big]}}{dt} & =-\left(p+\frac{q}{\lambda}\right)\overline{\big[S^\infty_{B_{1}}\big]}+\frac{q}{\lambda}\sum_{d=1}^{\infty}\frac{e^{-\lambda} \lambda^{d-1}}{(d-1)!}\left(\frac1{e^{-pt}(1-I^0)}\right)^{d-3}\left(\overline{\big[S^\infty_{B_{1}}\big]}\right)^{d-1}
 		\\ & = -\left(p+\frac{q}{\lambda}\right)\overline{\big[S^\infty_{B_{1}}\big]}+\frac{q}{\lambda}(1-I^0)^{2}e^{-2pt-\lambda}\sum_{d=0}^{\infty}\frac{\lambda^{d}}{d!}\left(\frac1{e^{-pt}(1-I^0)}\right)^{d}\left(\overline{\big[S^\infty_{B_{1}}\big]}\right)^{d}
 		\\ & = -\left(p+\frac{q}{\lambda}\right)\overline{\big[S^\infty_{B_{1}}\big]}+\frac{q}{\lambda}(1-I^0)^{2}e^{-2pt-\lambda+
 			\frac{\lambda}{ e^{-pt} (1-I^0)}\overline{\big[S^\infty_{B_{1}}\big]}}.
 	\end{aligned}
 	\label{eq:last}
 \end{equation}
 Let $y$ be given by~\eqref{eq:y=e^(pt)(1-beta_0)[S^infty_B_1]}.  
 %$y(t):=\frac1{e^{-pt}(1-I^0)}\overline{\big[S^\infty_{B_{1}}\big]}$.
 Then, using~(\ref{eq:last}),
 \[
 \begin{aligned}
 	\frac{dy}{dt}
 	& =py+\frac1{e^{-pt}(1-I^0)}\frac{d \overline{\big[S^\infty_{B_{1}}\big]}}{dt}
 	\\
 	& = py+\frac1{e^{-pt}(1-I^0)} \bigg( -\left(p+\frac{q}{\lambda}\right)\overline{\big[S^\infty_{B_{1}}\big]}+\frac{q}{\lambda}(1-I^0)^{2}e^{-2pt-\lambda+\frac{\lambda}{ e^{-pt} (1-I^0)}\overline{\big[S^\infty_{B_{1}}\big]}}\bigg)
 	\\ &  = 
 	py  -\left(p+\frac{q}{\lambda}\right)y+\frac{q}{\lambda}(1-I^0) e^{-pt-\lambda+ \lambda y}.
 \end{aligned}
 \]
 Therefore, $y$ satisfies~(\ref{eq:y-ER}), as claimed. Finally, by~\eqref{eq:pois}, \eqref{eq:s_by_dist}, and~\eqref{eq:s_k-ER}, 
 \begin{equation*}
 	\begin{aligned}
 		\overline{[S^\infty]}
 		&=\sum_{d=0}^{\infty} \mathbb{P}_{\rm deg}(d) \, \overline{[S^\infty_{B_{d}}]}
 		= \sum_{d=0}^{\infty}\frac{\lambda^{d}}{d!}e^{-\lambda}  \left(\frac1{e^{-pt}(1-I^0)}\right)^{d-1}
 		\left(\overline{\big[S^\infty_{B_{1}}\big]}\right)^{d}
 		\\ & =
 		(1-I^0) e^{-\lambda  -pt} \sum_{d=0}^{\infty}\frac{\lambda^{d}}{d!}  \left(\frac1{e^{-pt}(1-I^0)}\right)^{d}
 		\left(\overline{\big[S^\infty_{B_{1}}\big]}\right)^{d}  
 		%	  \\& 
 		%\\& 
 		=
 		(1-I^0) e^{-\lambda  -pt} \sum_{d=0}^{\infty}\frac{\lambda^{d}}{d!} y^{d} 
 		\\&	= 
 		(1-I^0) e^{-\lambda  -pt }e^{\lambda y}. 
 	\end{aligned}
 \end{equation*}
 Therefore, $f^{\rm ER} = 1-(1-I^0) e^{-\lambda  -pt+ \lambda y }$, as claimed. 
\end{proof}

With this, the proof of Theorem~\ref{thm:f^ER} is concluded.

\section{Explicit calculation of~$f^{ \text{\rm d-reg}}$}
\label{sec:f^d-regular-proof}

In this section we provide a proof of Theorem~\ref{thm:f^d-regular}.

\subsection{Number of cycles}

Let~$c_L(j; G^{\text{\rm d-reg}})$ denote the number of cycles with $L$~nodes that go through the node~$j$ in a graph~$G^{\text{\rm d-reg}}$.
The expected number of cycles of size~$L$ that go through the node~$j$
in~${\cal G}_{M}^{\text{\rm d-reg}}$  is
\begin{equation}
	C_L^{ \text{\rm d-reg}} :=
	\mathbb{E}_{{\cal G}_{M}^{\text{\rm d-reg}}} \big[ c_L(j) \big].
	%\frac{1}{M}\sum_{m =1}^M c_L(j; G^{\text{\rm d-reg}}).
	\label{eq:C_L-regular}
\end{equation} 
Since the family~${\cal G}_{M}^{\text{\rm d-reg}}$ is invariant under permutations of the nodes, $C_L^{ \text{\rm d-reg}}$ is independent of $j$.  %Therefore, $C_L^{ \text{\rm d-reg}}$ is also the expected number of cycles with $L$~nodes that go through a node in~${\cal G}_{M}^{\text{\rm d-reg}}$.
%	The expectation of~$c_L^{ \text{\rm d-reg}}$ in ${\cal G}_{M}^{\text{\rm d-reg}}$ has the following upper bound:
\begin{lemma}
	\label{lem:upper-bound-c_L^d-regular}
	Let $d \ge 3$. Then there exists a constant $\beta_d>0$, that depends only on~$d$, such that  
	%	the expectation of~$c_L^{ \text{\rm d-reg}}$ in ${\cal G}_{M}^{\text{\rm d-reg}}$ has the upper bound 
	\begin{equation}
		\label{eq:E_G[c_L]-d-regular-bound}
		C_L^{ \text{\rm d-reg}}
		< \frac{\kappa_{L}^{\text{\rm d-reg}}}{M} ,
		\qquad \kappa_{L}^{\text{\rm d-reg}} := \frac{1}{2}(\beta_d)^L, 
		\qquad L=3, \dots, M.
	\end{equation}
	%where $\kappa_{L}^{\text{\rm d-reg}} = \frac{1}{2}(\beta_d)^L$. 
\end{lemma}
\begin{proof} See Appendix~\ref{app:upper-bound-c_L^d-regular}.
\end{proof}

In the case of {\em short cycles}, one can replace the upper bound~\eqref{eq:E_G[c_L]-d-regular-bound} with 
an explicit limit: 
\begin{lemma}[\cite{Bollobas-80,Wormald-81}]
	\label{lem:E_G[c_L]-d-reg}
	%in ${\cal N}^{\rm ER} \in G\left(M,\frac{\lambda}{M}\right)$. 
	%Let $d \ge 3$. Then 
	For any fixed $L  \ge 3$,   
	%	The expectation of~$c_L^{ \text{\rm d-reg}}$ in ${\cal G}_{M,d}^{\rm regular}$ satisfies
	%over all $d$-regular graphs of the number of cycles of length~$L$ 
	%	that go through any given node satisfies  
	%in ${\cal N}^{\rm ER} \in G\left(M,\frac{\lambda}{M}\right)$. 
	\begin{equation}
		\label{eq:E_G[c_L]-d-regular}
		C_L^{ \text{\rm d-reg}} \sim 
		\frac{(d-1)^L}{2} \frac{1}{M}, \qquad M \to \infty.
	\end{equation}
\end{lemma}
%This asymptotic limit
%%~\eqref{eq:E_G[c_L^d]} for  $\mathbb{E}_{{\cal G}^{\rm ER}(M,\frac{\lambda}{M})}[c_L^d]$ 
%is confirmed 
%numerically in Figure~\ref{fig:circles_regular}.

%\begin{figure}[ht!]
%	\begin{center}
	%		\scalebox{0.6}{\includegraphics{Figures/ER/circles_reg.eps}}
	%		%		\scalebox{0.5}{\includegraphics{Figures/ER/Figure_26.eps}}
	%	\end{center}
%	\caption{The expected number of cycles of length~$L$ through a node  degree~$d$ in $d$-regular networks, as a function of the network size (circles). The solid line in the asymptotic limit~\eqref{eq:E_G[c_L]-d-regular}. A)~$d=4$, $L=3$. B)~$d=3$, $L=4$. C)~$d=4$, $L=4$. 
	%		D)~$d=5$, $L=4$. }
%	\label{fig:circles_regular}
%\end{figure}

%We also derive a lower bound for
%the average number of cycles through a node:
%\begin{lemma}
%	\label{lem:c^d-regular>3(d/2-1)}
%Let  $c^{ \text{\rm d-reg}}:=\sum_{L=3}^M c_L^{ \text{\rm d-reg}}$ denote 
%the average number of cycles through a node in~$ G^{\rm regular}$. 
%	Let $d \ge 3$. Then for any $d$-regular graph~$ G^{\rm regular}$, 
%	$
%	c^{ \text{\rm d-reg}}( G^{\rm regular})  \ge 3 \left(\frac{d}{2}-1 \right).
%	$
%\end{lemma}
%\begin{proof}  See Appendix~\ref{app:c^d-regular>3(d/2-1)}.
%\end{proof}

Lemma~\ref{lem:E_G[c_L]-d-reg} shows that on large $d$-regular random graphs, 
most of the nodes do not lie on short cycles. Indeed, 
{\em $d$-regular graphs tend to look like trees locally}, 
up to a $\log M$ distance from most nodes.

\subsection{Vanishing effect of cycles}

  	We can derive a simpler expression for~$[S^{ \text{\rm d-reg}}]$: 
\begin{lemma}
	Let 
	$[S^{ \text{\rm d-reg}}]  := 1- f^{ \text{\rm d-reg}}$. Then for any $ j \in \cal M$, 
	\begin{equation}
		\label{eq:S^d-regular-def}  
		[S^{ \text{\rm d-reg}}] = 
		\lim_{M \to \infty} \mathbb{E}_{{\cal G}_{M}^{\text{\rm d-reg}}} \big[\,[S_j]\, \big],  \qquad [S_j] = 1-f_j.
	\end{equation}
\end{lemma}
\begin{proof} 
	By~\eqref{eq:number_to_fraction-general} and~\eqref{eq:f^d-regular-def},
	\begin{equation*}
		%\label{eq:S^d-regular-def}  
		[S^{ \text{\rm d-reg}}] = 
		\lim_{M \to \infty} \mathbb{E}_{{\cal G}_{M}^{\text{\rm d-reg}}} \bigg[\frac{1}{M}\sum_{m =1}^M [S_m] \bigg]
		= 
		\lim_{M \to \infty} \frac{1}{M}\sum_{m =1}^M \mathbb{E}_{{\cal G}_{M}^{\text{\rm d-reg}}} \Big[[S_m] \Big]. 
	\end{equation*}
	Since the family~${\cal G}_{M}^{\text{\rm d-reg}}$ is invariant under  permutations of the nodes, $\mathbb{E}_{{\cal G}_{M}^{\text{\rm d-reg}}}[S_m]$ is independent of~$m$. Hence, the result follows.
\end{proof}

%	As in the case of ER networks, the effects of cycles vanish as $M \to \infty$. Therefore, in what follows,  we carry out the  calculations 
%	assuming that there are no cycles in the network, i.e., we use the funnel equality rather than the funnel inequality.
Consider the  $d$-regular network~${\cal N}^{ \text{\rm d-reg}}(G^{ \text{\rm d-reg}})$, where $G^{ \text{\rm d-reg}} \in {\cal G}_{M}^{\text{\rm d-reg}}$. Let $j \in \cal M$, 	
and denote by $\{k_{i}(j)\}_{i=1}^d$ the $d$~nodes connected to~$j$.
Define the networks~${\cal N}^{k_i,p_j}$  	for $i=1, \dots, d$ and the network~${\cal N}^{p_j}$ as in Theorem~\ref{thm:funnel-ER}. 
%Let 
%$    	[S^{ \text{\rm d-reg}}]  := 1- f^{ \text{\rm d-reg}}.
%$	
%Then
%\begin{equation}
%	\label{eq:S^d-regular-def}  
%	[S^{ \text{\rm d-reg}}] = 
%	\lim_{M \to \infty} \mathbb{E}_{{\cal G}_{M,d}^{\rm regular}} \left[\frac{1}{M}\sum_{j =1}^M [S_j] \right], \qquad [S_j] = 1-f_j.  
%\end{equation} 	
The following lemma shows that the effect of cycles become negligible as $M \to \infty$, so that one can effectively 
compute $[S^{ \text{\rm d-reg}}]$ using the funnel equality:   
\begin{lemma}
	\label{lem:S-funnel-equality-d-regular}
	\begin{equation}
		\label{eq:S-funnel-equality-d-regular}
		[S^{ \text{\rm d-reg}}]  = \lim_{M \to \infty} \mathbb{E}_{{\cal G}_{M,d}^{\rm regular}} \left[ \frac{1}{M}\sum_{j =1}^M  \frac{\prod_{i=1}^d [S^{k_i,p_j}_j]}{([S^{p_j}_j])^{d-1}}  \right],
	\end{equation}
	where
	\begin{equation}
		\label{eq:[S^pj_j](t;N^d-regular)=(1-beta0)e^-pt}
		[S^{p_j}_j]  \equiv  (1-I^0)e^{-pt}, \qquad 
		j\in {\cal M}.
	\end{equation}
\end{lemma}
\begin{proof}
Let $j \in\cal M$ and let 
$\tilde{q}:=\frac{q}{d}$.  
Fix $t>0$, and let $L^0 = L^0(t):= \left \lceil 2\tilde{q}t-1 \right \rceil$.
By Theorem~\ref{thm:funnel-ER}, 
and since $[S^{p_j}_j] \le 1$, see~\eqref{eq:[S^pj_j](t;N^ER)=(1-beta0)e^-pt}, 
$$
\begin{aligned}
	0 \le 
	[S_j] -	\frac{\prod_{i=1}^d [S^{k_i,p_j}_j]}{([S^{p_j}_j])^{d-1}}  
	%\\  &\
	\le
	2 
	\sum_{L=3}^{L^0(t)}c_L(j)
	\bigg(\frac{\tilde{q}}{p+\tilde{q}} \bigg)^{\left \lfloor{\frac{L+1}{2}}\right \rfloor} 
	+ 	2\sum_{L=L^0(t)+1}^{M}c_L(j) e^{-\left(p+\tilde{q}\right)t}
	\bigg(\frac{e\tilde{q}t}{\left \lfloor{\frac{L+1}{2}}\right \rfloor}\bigg)^{\left \lfloor{\frac{L+1}{2}}\right \rfloor},
\end{aligned}
$$
where $c_L(j)$ is the number of cycles of length~$L$ that go through~$j$.
%
%	
%	Fix $t>0$, and let $L^0 = L^0(t):= \left \lceil 2\tilde{q}t-1 \right \rceil$.
%	Then 
%	by~\eqref{eq:[S^p_j_j]<=1-regular}, \eqref{eq:funnel-ineq-regular-d}, and~\eqref{eqs:funnel-ineq-regular-d-upper-bound},
%	we have that  
%	$$
%	\begin{aligned}
%		0 \le 
%		[S_j] -	\frac{\prod_{i=1}^d [S^{k_i,p_j}_j]}{([S^{p_j}_j])^{d-1}}  
%		%\\  &\
%		\le
%		2 
%		\sum_{L=3}^{L^0(t)}c_L(j)
%		\left(\frac{\tilde q}{p+\tilde q} \right)^{\left \lfloor{\frac{L+1}{2}}\right \rfloor} 
%		+ 	2\sum_{L=L^0(t)+1}^{M}c_L(j) e^{-(p+\tilde{q})t}
%		\left(\frac{e\tilde{q}t}{\left \lfloor{\frac{L+1}{2}}\right \rfloor}\right)^{\left \lfloor{\frac{L+1}{2}}\right \rfloor},
%	\end{aligned}
%	$$
%	where $c_L(j)$ is the number of cycles of length~$L$ that go through~$j$. 
	Averaging these inequalities %over all $j \in \cal M$ and then averaging 
	over ${\cal G}_{M,d}^{\rm regular}$ gives, see~\eqref{eq:C_L-regular}, 
	\begin{subequations}
		\label{eqs:E_G[S]}
		\begin{equation}
	\begin{aligned}
		0 
		\le 
		\mathbb{E}_{{\cal G}_{M}^{\text{\rm d-reg}}} \bigg[ \, %\frac{1}{M}\sum_{j =1}^M 
		[S_j] -	%\frac{1}{M}\sum_{j =1}^M 
		\frac{\prod_{i=1}^d [S^{A_i,p_j}_j]}{([S^{p_j}_j])^{d-1}}    \, \bigg]
		%\\  &
		\le
		2  \sum_{L=3}^{L_0(t)} 
		C_L^{\text{\rm d-reg}}
		\left(\frac{\tilde q}{p+\tilde q} \right)^{\left \lfloor{\frac{L+2}{2}}\right \rfloor} 
		+ 
		2 \!\!\!\!\!\! \sum_{L=L_0(t)+1}^{M} 
		C_L^{\text{\rm d-reg}}
		%\big[c_L^{ \text{\rm d-reg}}\big] %e^{-(p+\tilde{q})t}
		\left(\frac{e\tilde{q}t}{\left \lfloor{\frac{L+2}{2}}\right \rfloor}\right)^{\left \lfloor{\frac{L+2}{2}}\right \rfloor}.
	\end{aligned}
%	\vspace{2mm}
\end{equation}
	%	where $c_L^{ \text{\rm d-reg}}$ is given by~\eqref{eq:c_L-regular}. 
		%	where  $c_L(d, G)$ is the number of different cycles with $L$~nodes that go through a node of degree~$d$.
		%		Recall that, see~\eqref{eq:E_G[c_L]-d-regular}, 
		%		\begin{equation*}
			%			\mathbb{E}_{{\cal G}_{M,d}^{\rm regular}} \big[c_L^{ \text{\rm d-reg}}\big] \sim \frac{ d (d-1)}{2} \lambda^{L-2} \frac1{M}, \qquad M \to \infty.
			%		\end{equation*}
	Using the upper bound~\eqref{eq:E_G[c_L]-d-regular-bound} 
	for~$C_L^{ \text{\rm d-reg}}$, we have
		\begin{equation}
	\begin{aligned}
		\sum_{L=3}^{L_0(t)} C_L^{\text{\rm d-reg}}
		\left(\frac{\tilde q}{p+\tilde q} \right)^{\left \lfloor{\frac{L+2}{2}}\right \rfloor}
		<
		\sum_{L=3}^{L_0(t)}\frac{\kappa_{L}^{\text{\rm d-reg}}}{M}
		\left(\frac{\tilde q}{p+\tilde q} \right)^{\left \lfloor{\frac{L+2}{2}}\right \rfloor}
		%\\ &   
		= O\left(\frac1M \right), \qquad M \to \infty,
	\end{aligned}
\end{equation}
and  
\begin{equation}
	\begin{aligned}
		&	\sum_{L=L_0(t)+1}^{M} C_L^{\text{\rm d-reg}} 
		\left(\frac{e\tilde{q}t}{\left \lfloor{\frac{L+2}{2}}\right \rfloor}\right)^{\left \lfloor{\frac{L+2}{2}}\right \rfloor}
		< 
		\sum_{L=L_0(t)+1}^{M} \frac{\kappa_{L}^{\text{\rm d-reg}}}{M}
		\left(\frac{e\tilde{q}t}{\left \lfloor{\frac{L+2}{2}}\right \rfloor}\right)^{\left \lfloor{\frac{L+2}{2}}\right \rfloor}
		\\ & \qquad\qquad\qquad < 
		\frac1{M}  \sum_{L=L_0(t)+1}^{\infty} \frac{(\beta_d)^L}{2}
		\left(\frac{e\tilde{q}t}{\left \lfloor{\frac{L+2}{2}}\right \rfloor}\right)^{\left \lfloor{\frac{L+2}{2}}\right \rfloor}
		= O\left(\frac1M \right), \qquad M \to \infty.
	\end{aligned}
\end{equation}
	\end{subequations}
	Letting $M \to \infty$   in relations~\eqref{eqs:E_G[S]} and using~\eqref{eq:S^d-regular-def} proves~\eqref{eq:S-funnel-equality-d-regular}.
 \end{proof}

%\begin{equation*}
%%	\label{eq:lim_M->infty_E_G[S_B_d]}
%	\overline{\left[S^\infty_{B_{d}}\right]}:= 
%		\lim_{M \to \infty} \mathbb{E}_{G\left(M,\frac{\lambda}{M}\right)} \left[\, [S_{B_{d}}] \, \right] = 		\lim_{M \to \infty}
%		 \mathbb{E}_{G\left(M,\frac{\lambda}{M}\right)} \left[\, 
%		\frac{1}{|B_{d}|}\sum_{j\in B_{d}} \frac{\prod_{i=1}^d [S^{k_i,p_j}_j]}{([S^{p_j}_j])^{d-1}}    \right] ,
%\end{equation*}
%as needed.

Recall that most nodes belong, on average, to at least 
one cycle of size $\log M$, and that a strict funnel inequality holds for any node that lies on a cycle
(Theorem~\ref{thm:funnel-ER}). Nevertheless,
Lemma~\ref{lem:S-funnel-equality-d-regular} shows that 
one can compute $[S_j]$ using the funnel {\em equality}, see~\eqref{eq:S^d-regular-def}. This is because as 
the cycle length~$L$ increases, the number of cycles increases exponentially in~$L$, see~\eqref{eq:E_G[c_L]-d-regular}, but 
the effect of cycles on the adoption probability decays 
at a super-exponential ${L^{-L}}$~rate,
see~\eqref{eq:funnel-upper-bound-Ncycle}.
Therefore, the overall effect of cycles vanishes as $M \to \infty$.

\subsection{Solving the master equations}

Next, we simplify the right-hand side of~\eqref{eq:S-funnel-equality-d-regular}:
\begin{lemma}
	\label{lem:end-pf-d-regular}
	For any $j \in \cal M$, 
	\begin{equation}
		\label{eq:end-pf-d-regular}
		\lim_{M \to \infty} \mathbb{E}_{{\cal G}_{M,d}^{\rm regular}} \left[  \frac{\prod_{i=1}^d [S^{k_i,p_j}_j]}{([S^{p_j}_j])^{d-1}}  \right]
		=
		\left(\frac{1}{e^{-pt}(1-I^0)}\right)^{d-1}
		\prod_{i=1}^d 
		\lim_{M \to \infty} 
		\mathbb{E}_{{\cal G}_{M,d}^{\rm regular}} \left[ [S^{k_i,p_j}_{j}]  \right].
	\end{equation}
\end{lemma}
\begin{proof}
	By~\eqref{eq:[S^pj_j](t;N^d-regular)=(1-beta0)e^-pt}, 
	\begin{subequations}
		\label{eqs:end-pf-d-regular}
		\begin{equation}
			\mathbb{E}_{{\cal G}_{M,d}^{\rm regular}} \left[  \frac{\prod_{i=1}^d [S^{k_i,p_j}_j]}{([S^{p_j}_j])^{d-1}}  \right]
			= 
			\left(\frac{1}{e^{-pt}(1-I^0)}\right)^{d-1}	
			\mathbb{E}_{{\cal G}_{M,d}^{\rm regular}} 
			\left[ \prod_{i=1}^d [S^{k_i,p_j}_j]  \right].
		\end{equation}
		
%		In addition, since $\mathbb{E}_{{\cal G}_{M,d}^{\rm regular}} \left[\prod_{i=1}^d  [S^{k_i,p_j}_j] \right] $ is independent of~$j$,
%		it follows that  for any $j_0 \in \cal M$, 
%		\begin{equation}
%			\frac{1}{M}\sum_{j =1}^M \mathbb{E}_{{\cal G}_{M,d}^{\rm regular}} \left[ \prod_{i=1}^d 	 [S^{k_i,p_j}_j] \right]
%			= 
%			\mathbb{E}_{{\cal G}_{M,d}^{\rm regular}} \left[ \prod_{i=1}^d  [S^{k_i,p_j}_j] \right] .
%		\end{equation}
The same proof as in Lemma~\ref{lem:lemma3} shows that 
%	\begin{equation}
%%	\label{eq:end-pf-d-regular}
%	\lim_{M \to \infty} \mathbb{E}_{{\cal G}_{M,d}^{\rm regular}} \left[  \frac{\prod_{i=1}^d [S^{k_i,p_j}_j]}{([S^{p_j}_j])^{d-1}}  \right]
%	=
%	\left(\frac{1}{e^{-pt}(1-I^0)}\right)^{d-1}
%	\prod_{i=1}^d 
%	\lim_{M \to \infty} 
%	\mathbb{E}_{{\cal G}_{M,d}^{\rm regular}} \left[ [S^{k_i,p_j}_{j}]  \right].
%\end{equation}
%
%
%		Since $[S^{k_i,p_j}_j]$ is completely determined by~$[S_{k_i}]$,
%		see~\eqref{eq:d[S^(k_i,p_j)_j]_dt},
%		and since~$\{[S_{k_{i}}] \}_{i=1}^d$
%		become independent as $M \to \infty$, 
%		the probabilities $\left\{[S^{k_i,p_j}_j]\right\}_{i=1}^d$
%		become independent in~$\lim_{M \to \infty} {\cal G}_{M,d}^{\rm regular}$ as well.  
%		Hence,
		\begin{equation}
			\label{eqs:end-pf-d-regular-product}
			\lim_{M \to \infty} 
			%	 \frac{1}{M}\sum_{j =1}^M 
			\mathbb{E}_{{\cal G}_{M,d}^{\rm regular}} \left[ \prod_{i=1}^d [S^{k_i,p_j}_j]  \right]
			=
			%	 \frac{1}{M}\sum_{j =1}^M 
			\prod_{i=1}^d 
			\lim_{M \to \infty} 
			\mathbb{E}_{{\cal G}_{M,d}^{\rm regular}} \left[ [S^{k_i,p_j}_j]  \right].
		\end{equation}
	\end{subequations}
	Combining relations~\eqref{eqs:end-pf-d-regular} yields the result. 
 \end{proof}

From equations~\eqref{eq:S-funnel-equality-d-regular} and~\eqref{eq:end-pf-d-regular}, we have that
%\begin{subequations}
%	\label{eqs:S-funnel-equality-d-regular-interim}
\begin{equation}
	\label{eq:S-funnel-equality-d-regular-interim}
	[S^{ \text{\rm d-reg}}]  = 
	\left(\frac{1}{e^{-pt}(1-I^0)}\right)^{d-1}
	\prod_{i=1}^d 
	\lim_{M \to \infty} 
	\mathbb{E}_{{\cal G}_{M,d}^{\rm regular}} \left[ [S^{k_i,p_j}_j]  \right].
\end{equation}
Let $y(t)$ denote the expected nonadoption probability of a node in infinite $d$-regular networks, after we delete $d-1$ of its edges.
Then
\begin{equation}
	\label{eq:y-def-d-regular-pf}
	y(t) =
	\lim_{M \to \infty}\mathbb{E}_{{\cal G}_{M,d}^{\rm regular}} \left[ [{S_j^{A_{i}, p_j}}]  \right], 
	\qquad i\in \{1, \dots, d\}, \quad j_0 \in \cal M.  
\end{equation}
By~\eqref{eq:S-funnel-equality-d-regular-interim} and~\eqref{eq:y-def-d-regular-pf},
\begin{equation}
	\label{eq:S^d-regular=y^d}
	[S^{ \text{\rm d-reg}}] = 	\left(\frac{1}{e^{-pt}(1-I^0)}\right)^{d-1} y^d.
\end{equation}
The auxiliary function $y(t)$ satisfies the following ODE:
\begin{lemma}
	\begin{equation}  
		\label{eq:dy_dt_f-regular}
		\frac{dy}{dt} +\left(p+\frac{q}{d}\right)y = \frac{q}{d} 
		\left(\frac{1}{e^{-pt}(1-I^0)}\right)^{d-3} y^{d-1}.
	\end{equation}
\end{lemma}
\begin{proof}
	%	Finally, since~${\cal N}^{k_i,p_j}$  is a $d$-regular network in ${\cal G}_{M,d}^{\rm regular}$ in which $j$~is a degree-one node,
	%	%	 and it can only be influenced by nodes in~$A_i$. Therefore,
	%	%	 $j$~is a degree-one node in an  ER network in $G\left(\frac{M}{d},\alpha\right)$
	%	\begin{equation}
		%		\label{eq:E_G[[S^A_i,p_j_j]=[S_B_1]}
		%		\lim_{M \to \infty} 
		%		\mathbb{E}_{G\left(M,\frac{\lambda}{M}\right)}\left[\, 
		%		\frac{1}{|B_{d}|}\sum_{j\in B_{d}} [S^{k_i,p_j}_j]\right] = \overline{\big[S^\infty_{B_{1}}\big]}, \qquad i=1, \dots, d.
		%	\end{equation} 
	%\end{subequations}
	%%	  Hence, using~\eqref{eq:E_G[[S^A_i,p_j_j]=[S_B_1]}, 
	%%	  \begin{equation}
		%	%	  	\label{eq:lim_M->infty_E_G[product]}
		%	%	  \lim_{M\to \infty}	\mathbb{E}_G\left[\prod_{i=1}^d [S^{k_i,p_j}_j] \right] = 
		%	%	  	\prod_{i=1}^d \lim_{M\to \infty} \mathbb{E}_G\left[ [S^{k_i,p_j}_j] \right]
		%	%	  	 =  \left(\overline{\big[S^\infty_{B_{1}}\big]}\right)^{d}.
		%	%	  \end{equation}
	%The result follows from~\eqref{eq:lim_M->infty_E_G[S_B_d]}
	%and~\eqref{eqs:end-pf-ER}.
	%	By the funnel equality, see~\eqref{eq:funnel_equality-A-p-K>2}, 
	% \begin{subequations}
		%	\label{eqs:y(t)-d-regular}
		%	\begin{equation}
			%[S_j] = (1-I^0)^{-(d-1)}e^{(d-1)pt} \prod_{i=1}^d [S_j^{A_{i}, p_j}],
			%\end{equation}	
			%where $[S_j^{A_{i}, p_j}]$ is the nonadoption probability of~$j$ in 
			%network~${\cal N}^{A_{i}, p_j}$ in which we delete $d-1$ edges of~$j$
			%and leave only its edge to~$k_i$.		
			%\,\footnote{By symmetry, 
				%$
				%[\widetilde{S_j}] := [S_j^{A_{i}, p_j}]
				%$
				%is independent of~$i$.} 
			By the indifference principle, we can compute 
			$[S_j^{A_{i}, p_j}]$ on an equivalent network, denoted by~$\widetilde{{\cal N}}$,
			in which we delete the directional edge from $j$ to $k_i$ (but leave the directional edge from~$k_i$ to $j$), i.e.,
			\begin{equation}
				[S_j^{A_{i}, p_j}] = [\widetilde{S_j^{A_{i}, p_j}}],
			\end{equation}
			where $[\widetilde{\, \cdot \,} ]$ denote nonadoption probabilities in~$\widetilde{{\cal N}}$. 
			In network~$\widetilde{{\cal N}}$, $j$ is a degree-one node
			which is directly influenced by~$k_i$ but is not 
			influential to $k_i$. Therefore, 
			by~\cite[Theorem~4.1]{Bass-Stochastic-IC-22},  the master equation for~$[\widetilde{S_j^{A_{i}, p_j}}]$ is given by
			\begin{equation}
				\label{eq:d_dt_[widetilde_S_j^A_i,p_j]}
				\frac{d}{dt} [\widetilde{S_j^{A_{i}, p_j}}]+\left(p+\frac{q}{d}\right)[\widetilde{S_j^{A_{i}, p_j}}] 
				=
				\frac{q}{d} (1-I^0) e^{-pt } \, [\widetilde{S_{k_i}}],\quad t>0,
				\qquad [\widetilde{S_j^{A_{i}, p_j}}](0)=1-I^0.
			\end{equation}
			
			Let
			$$
			y_M(t):=
			\mathbb{E}_{{\cal G}_{M,d}^{\rm regular}} \left[ [\widetilde{S_j^{A_{i}, p_j}}]  \right].
			$$
			By~\eqref{eq:d_dt_[widetilde_S_j^A_i,p_j]}, 
			\begin{equation*}
				%	\label{eq:d_dt_[widetilde_S_j^A_i,p_j]}
				\frac{d}{dt} y_M
				+\left(p+\frac{q}{d}\right)y_M = \frac{q}{d} (1-I^0) e^{-pt } \, 
				\mathbb{E}_{{\cal G}_{M,d}^{\rm regular}} \left[ [\widetilde{S_{k_i}}]\right],\quad t>0,
				\qquad y_M(0)=1-I^0.
			\end{equation*}
			Letting $M \to \infty$ and using~\eqref{eq:y-def-d-regular-pf} gives
			\begin{subequations}
				\label{eqs:S-funnel-equality-d-regular-interim}
				\begin{equation}
					%	\label{eq:d_dt_[widetilde_S_j^A_i,p_j]}
					\frac{d}{dt} y
					+\left(p+\frac{q}{d}\right)y =
					\frac{q}{d} (1-I^0) e^{-pt } \, 
					\lim_{M \to \infty} \mathbb{E}_{{\cal G}_{M,d}^{\rm regular}} \left[ [\widetilde{S_{k_i}}]\right],\quad t>0,
					\qquad y(0)=1-I^0.
				\end{equation}

				In the network~$\widetilde{{\cal N}}$, $k_i$ is a node of degree $d-1$ in an otherwise $d$-regular network. 
				Denote by $\{m_{n}\}_{n=1}^{d-1}$ the $d-1$ nodes connected to~$k_i$. 
%				Let $\{B_n\}_{n=1}^{d-1}$ be a partition of ${\cal M} \setminus \{j_0,k_i\}$, such that $m_{n} \in  B_n$ for $n=1, \dots, d-1$.  
%				
Using the same steps as in the proof of Lemma~\ref{lem:S-funnel-equality-d-regular}, 
				we have that 
				\begin{equation}
					%	\label{eq:S-funnel-equality-d-regular}
					\lim_{M \to \infty} \mathbb{E}_{{\cal G}_{M,d}^{\rm regular}} \left[ [\widetilde{S_{k_i}}] \right]
					= 
					\lim_{M \to \infty} \mathbb{E}_{{\cal G}_{M,d}^{\rm regular}} \Bigg[% \frac{1}{M}\sum_{j =1}^M  
					\frac{\prod_{n=1}^{d-1} \Big [\widetilde{S_{k_i}^{m_n, p_{k_i}}} \Big]}{\Big([\widetilde{S^{p_{k_i}}_{k_i}}]\Big)^{d-2}}  \Bigg],
				\end{equation}
				where $\Big [\widetilde{S_{k_i}^{m_n, p_{k_i}}} \Big]$ is the nonadoption probability of~$k_i$ in the network~$\widetilde{{\cal N}^{m_n, p_{k_i}}}$, in which we delete $d-2$ edges of~$k_i$
				and leave only its edge to~$m_{n}$.
				The same arguments as in the proof of Lemma~\ref{lem:end-pf-d-regular}
				show that
				\begin{equation}
					%	\label{eq:S-funnel-equality-d-regular}
					%	\lim_{M \to \infty} \mathbb{E}_{{\cal G}_{M,d}^{\rm regular}} \left[ [\widetilde{S_{k_i(j_0)}}] \right] 
					\lim_{M \to \infty} \mathbb{E}_{{\cal G}_{M,d}^{\rm regular}} \Bigg[% \frac{1}{M}\sum_{j =1}^M  
					\frac{\prod_{n=1}^{d-1} \Big [\widetilde{S_{k_i}^{m_n, p_{k_i}}} \Big]}{\Big([\widetilde{S^{p_{k_i}}_{k_i}}]\Big)^{d-2}}  \Bigg]
					= \left(\frac{1}{e^{-pt}(1-I^0)}\right)^{d-2} \prod_{n=1}^{d-1}
					\lim_{M \to \infty}
					\mathbb{E}_{{\cal G}_{M,d}^{\rm regular}} \left[ \Big [\widetilde{S_{k_i}^{m_n, p_{k_i}}} \Big] \right] .
				\end{equation}
				% funnel equality, see~\eqref{eq:funnel_equality-A-p-K>2}, {\bf correct!!!}
				%	\begin{equation}
					%[\widetilde{S_{k_i}}] = (1-I^0)^{-(d-2)} e^{(d-2)pt} \prod_{n=1}^{d-1} \Big [\widetilde{S_{k_i}^{B_n, p_{k_i}}} \Big],
					%\end{equation}	
					Since $\Big [\widetilde{S_{k_i}^{m_n, p_{k_i}}} \Big]$ is the nonadoption probability of a degree-one node in an otherwise infinite $d$-regular network, it follows that
					\begin{equation}
						y =
						%   	\lim_{M \to \infty}\mathbb{E}_{{\cal G}_{M,d}^{\rm regular}} \left[ [\widetilde{S_j^{A_{i}, p_j}}]  \right]  
						%   	 = 
						\lim_{M \to \infty} \mathbb{E}_{{\cal G}_{M,d}^{\rm regular}} \left[  \Big [\widetilde{S_{k_i}^{m_n, p_{k_i}}} \Big] \right],
						\qquad i=1, \dots, d, \quad n=1, \dots, d-1  .
					\end{equation}
				\end{subequations}
				The result follows from relations~\eqref{eqs:S-funnel-equality-d-regular-interim}.
			 \end{proof}

			Let $[S]:=[S^{ \text{\rm d-reg}}]$. Then using~\eqref{eq:S^d-regular=y^d} and~\eqref{eq:dy_dt_f-regular},
			$$
			\begin{aligned}
				\frac{d[S]}{dt} & = (d-1)p [S]+ d \frac{[S]}{y}\frac{dy}{dt} =
				(d-1)p [S]- d \frac{[S]}{y} \left( \left(p+\frac{q}{d}\right)y - \frac{q}{d}  \left(\frac{1}{e^{-pt}(1-I^0)}\right)^{d-3} \! y^{d-1}\right)	 
				\\ & = [S] \bigg( -p-q + q \left(\frac{1}{e^{-pt}(1-I^0)}\right)^{d-3}  y^{d-2} \bigg)
				\\ &= [S] \bigg( -p-q + q \left(\frac{1}{e^{-pt}(1-I^0)}\right)^{d-3}  \bigg([S]  \left(\frac{1}{e^{-pt}(1-I^0)}\right)^{-(d-1)}  \bigg)^{\frac{d-2}{d}} \bigg)  
				\\ & = [S] \bigg( -p-q + q\left(\frac{1}{e^{-pt}(1-I^0)}\right)^{-\frac{2}{d}} [S]^{\frac{d-2}{d}} \bigg).
			\end{aligned}
			$$
			This concludes the proof of Theorem~\ref{thm:f^d-regular}.

\appendix	
	
	\section{Proof of Lemma~\ref{lem:f^ER-increasing-in-lambda}}
	\label{app:f^ER-increasing-in-lambda}
	
	We can rewrite~\eqref{eqs:f^ER} as
	\begin{subequations}
		\label{eqs:f^ER-z(t)}
		\begin{equation}
			%	\lim_{M \to \infty}	\overline{\left[S\right]}(t)
			f^{\rm ER}(t;p,q,\lambda,I^0) = 1-
			(1-I^0)e^{-pt-z(t)}, \qquad t \ge 0,
			\label{eq:f^ER-z(t)}
		\end{equation}
		where $z: = \lambda (1-y)$ is the solution of the equation
		\begin{equation}
			\frac{dz}{dt} =-q\left(\frac{z}{\lambda}-1+(1-I^0)e^{-pt-z}\right), ~~ t \ge 0,
			\qquad  z(0)=0.
			\label{eq:z-ER}
		\end{equation}
	\end{subequations}
%	Since the right-hand side of~\eqref{eq:z-ER} is increasing with $\lambda$, so is $z(t;\lambda)$.  To prove this rigorously, 
     We can prove that $z$ is monotonically increasing in~$\lambda$ for all times, as follows. 
	Substituting $t=0$ in~\eqref{eq:z-ER} gives 
	$\frac{dz}{dt}(0) = q I^0$, which is independent of~$\lambda$.
	Differentiating~\eqref{eq:z-ER}  and substituting $t=0$ gives
	$$
	\frac{d^2 z}{dt^2}(0)
	= 
	-\frac{q}{\lambda}
	\frac{dz}{dt} (0)+\{\text{terms independent of $\lambda$}\}
	=  -\frac{q^2 I^0}{\lambda}  +\{\text{terms independent of $\lambda$}\}.
	$$
	Since $\frac{d^2 z}{dt^2}(0)$ is monotonically 
	increasing in~$\lambda$,  it follows that 
	$z(t;\lambda)$ is monotonically increasing in~$\lambda$ for $0 <t \ll 1$.
	
	Assume by contradiction that $z$ is not monotonically increasing in~$\lambda$ for all $t>0$. Then there exist $\lambda_1<\lambda_2$ and $t_0>0$ such that 
	$z(t_0;\lambda_1) \ge  z(t_0;\lambda_2)$.
	Therefore, there exists $0<t_1 \le t_0$ such that 
	$z(t_1;\lambda_1)=z(t_1;\lambda_2)$
	and 
	$z(t;\lambda_1) < z(t;\lambda_2)$ for $0<t <t_1$.
	Hence,
	$\frac{d z}{dt} (t_1;\lambda_1) \ge \frac{dz}{dt}(t_1; \lambda_2)$.
	Since, however, $z(t_1;\lambda_1)=z(t_1;\lambda_2)>0$, it follows from~\eqref{eq:z-ER} that $\frac{d z}{dt} (t_1;\lambda_1) < \frac{dz}{dt}(t_1; \lambda_2)$, which is a contradiction. Therefore, we conclude that
	$z$ is monotonically increasing in~$\lambda$ for all times. 
		 Hence, by~\eqref{eq:f^ER-z(t)}, so is~$f^{\rm ER}(t;\lambda)$. 	
	
		\section{Proof of Lemma~\ref{lem:f^ER_Bass-bounds}}
	\label{app:f^ER_Bass-bounds}
	
		When $\lambda=0$,  all the network nodes are isolated, and so we have~\eqref{eq:f^ER_lambda=0}.
	%see~\eqref{eq:isolated-nodes-network-f}. 
	By~\eqref{eqs:f^ER}, the ODE for $f:=f^{\rm ER}$ is  
	\begin{equation}
		\label{eq:df^ER_Bass_dt-lambda}
		\begin{aligned}
			\frac{df}{dt} 
			& =  (1-I^0)e^{-pt-\lambda(1-y)} \Big(p-\lambda \frac{dy}{dt}\Big)
			= \left(1-f\right)  \left(p-q \left(-y+(1-I^0)e^{-pt-\lambda(1-y)}\right)\right)
			\\ &
			= \left(1-f\right) \Big(p+q \left(y-1+f\right)\Big), \qquad f(0)=I^0.
		\end{aligned}
	\end{equation}
	Since $y(t) \le 1$, see Lemmas~\ref{lem:y(t)-properties-ER} and~\ref{lem:y_SI-ER} below, it follows from~\eqref{eq:y-ER}
	that $\frac{dy}{dt}> -\frac{q}{\lambda} y> -\frac{q}{\lambda}$. Therefore, 
	$1 \ge y(t) \ge  1 -t\frac{q}{\lambda}$, and so 
	$	\lim_{\lambda \to \infty} y(t) \equiv 1$.
	Hence,  as $ \lambda \to \infty$, 
	equation~\eqref{eq:df^ER_Bass_dt-lambda} reduces to  equation~\eqref{eq:f^complete_infty-ODE} for~$f^{\rm compart}$.
	%Therefore, we proved~\eqref{eq:lim_lambda_to_infty_f^ER}. 
	Finally, since~$f^{\rm ER}$ is monotonically increasing in~$\lambda$
	(Lemma~\ref{lem:f^ER-increasing-in-lambda}),  the lower and upper  bounds~\eqref{eq:bounds-f^ER} follow.
	% from~\eqref{eq:f^ER_lambda=0} and~\eqref{eq:lim_lambda_to_infty_f^ER}.
%	
%	
%	When $\lambda=0$,  all the network nodes are isolated. Hence, 
%	\begin{equation}
%		\label{eq:f^ER_Bass-lambda=0}
%		f^{\rm ER}(t;p,q,I^0,\lambda=0) = 1-(1-I^0)e^{-pt},
%	\end{equation}
%	see~\eqref{eq:s_0-ER}. 
%	By~\eqref{eqs:f^ER}, the ODE for $f:=f^{\rm ER}$ is  
%	\begin{equation}
%		\label{eq:df^ER_Bass_dt-lambda}
%		\begin{aligned}
%			\frac{df}{dt} 
%			& =  (1-I^0)\Big(p-\lambda \frac{dy}{dt}\Big)
%			e^{-pt-\lambda(1-y)}
%			=  \left(p-q \left(-y+(1-I^0)e^{-pt-\lambda(1-y)}\right)\right) \left(1-f\right)
%			\\ &
%			=  \left(p+q \left(y-1+f\right)\right)\left(1-f\right), \qquad f(0)=I^0.
%		\end{aligned}
%	\end{equation}
%	If we let $ \lambda \to \infty$ in  equation~\eqref{eq:y-ER}, we get that 
%	$	\lim_{\lambda \to \infty} y(t) \equiv 1$. 
%	Therefore, as $ \lambda \to \infty$, equation~\eqref{eq:df^ER_Bass_dt-lambda} becomes 
%\begin{equation}
%	\label{eq:f^complete_infty-ODE}
%	\frac{df}{dt} = (p+qf)(1-f), \qquad f(0) = I^0.
%\end{equation}
%The solution of this equation is given by~$f^{\rm compart}$, see~\eqref{eq:f^compart}. 
%	Therefore, we conclude that
%	\begin{equation}
%		\label{eq:lim_lambda_to_infty_f^ER}
%		\lim_{\lambda \to \infty} f^{\rm ER}(t;p,q,I^0,\lambda) = f^{\rm compart}(t;p,q,I^0).
%	\end{equation}
%	Since $f^{\rm ER}$ is monotonically increasing in~$\lambda$
%	(Lemma~\ref{lem:f^ER-increasing-in-lambda}), the result follows from~\eqref{eq:f^ER_Bass-lambda=0} and~\eqref{eq:lim_lambda_to_infty_f^ER}.
%	% \end{proof}
	
	\section{Proof of Lemma~\ref{lem:f^ER_Bass-monotone}}
	\label{app:f^ER_Bass-monotone}
	
	We first show that the auxiliary function~$y_{\rm Bass}(t)$ is monotonically decreasing from $y_{\rm Bass}(0)=1$ to~$y_{\rm Bass}(\infty)=0$:
	\begin{lemma}
		\label{lem:y(t)-properties-ER}
		Let~$y_{\rm Bass}(t)$ be the solution of~\eqref{eq:y-ER-Bass}. Then 
		\begin{equation}
			\label{eq:y(t)-properties-ER}
			%	y(0) = 1, \qquad 
			\frac{dy_{\rm Bass}}{dt}<0, \quad t>0   \qquad \text{and} \qquad \lim_{t \to \infty} y_{\rm Bass} = 0.
		\end{equation}
		In particular, $0<y_{\rm Bass}(t)<1$ for $t>0$.
	\end{lemma}
	\begin{proof}
		%	See Section~\ref{sec:y(t)-properties-ER}.
		%\end{proof}
		%
		%
		%
		%
		%\section{Proof of Lemma~\ref{lem:y(t)-properties-ER}}
		%\label{sec:y(t)-properties-ER}
		%
		Substituting $t=0$ in~\eqref{eq:y-ER-Bass} gives, for $y = _{\rm Bass}$, 
		\begin{equation*}
			%\label{eq:dy_dt(0)=-q*beta_0}
			\frac{dy}{dt}(0)=-\frac{q}{\lambda} (1-1) = 0.
		\end{equation*}
		Let~$\tilde{t} \ge 0$ for which $\frac{dy}{dt}(\tilde{t})=0$.
		Differentiating~\eqref{eq:y-ER-Bass}
		and substituting $t=\tilde{t}$ gives
		\begin{equation}
			\label{eq:dy_dt^2(0)=-q0*(1-beta_0}
			\begin{aligned}
				\frac{d^{2}y}{dt^{2}}(\tilde{t})= & \frac{q}{\lambda}\left(-\frac{dy}{dt}(\tilde{t})+\Big(-p+\lambda\frac{dy}{dt}(\tilde{t})\Big)e^{-p\tilde{t}-\lambda(1-y(\tilde{t}))}\right)\\
				= & 
				-\frac{q}{\lambda}p e^{-p\tilde{t}-\lambda(1-y(\tilde{t}))} < 0.
			\end{aligned}
		\end{equation}
		In particular, $\frac{d^{2}y}{dt^{2}}(0)<0$. Therefore,  
		%	From~\eqref{eq:dy_dt(0)=-q*beta_0} and~\eqref{eq:dy_dt^2(0)=-q0*(1-beta_0} it follows that 
		%	under conditions~\eqref{eq:p,q,I^0>=0-ER} and~\eqref{eq:p,I^0>0-ER},   either $\frac{dy}{dt}(0)<0$, or $\frac{dy}{dt}(0)=0$
		%	and  $\frac{d^{2}y}{dt^{2}}(t_{0})<0$. 
		%	Therefore,  
		$\frac{dy}{dt}<0$ for $0<t \ll 1$. 
		
		Assume by contradiction that there exists~$t_{0} > 0$ for which 	\begin{subequations}
			\label{eqs:dy_dt>0<0-ER}
			\begin{equation}
				\label{eq:dy_dt=0-ER_Bass}
				\frac{dy}{dt}(t_{0})=0.
			\end{equation}
			If $\frac{dy}{dt}$ vanishes at more than one positive time, let~$t_0$ be 
			the smallest positive such time. Therefore, %for $0<\epsilon \ll 1$,
			%	\begin{subequations}
				%		\label{eqs:dy_dt>0<0-ER}
				\begin{equation}
					\frac{dy}{dt}<0, \qquad 0 < t<t_{0}.
				\end{equation}
			\end{subequations} 
			From relations~\eqref{eqs:dy_dt>0<0-ER} it follows that
			$
			\frac{d^{2}y}{dt^{2}}(t_0) \ge 0.
			$
			From~\eqref{eq:dy_dt^2(0)=-q0*(1-beta_0} and~\eqref{eq:dy_dt=0-ER_Bass}, however, we have that
			$
			\frac{d^{2}y}{dt^{2}}(t_0) < 0,
			$
			%	In addition, from~\eqref{eq:dy_dt(0)=-q*beta_0} 
			%		and~\eqref{eq:dy_dt^2(0)=-q0*(1-beta_0} it follows that 
			%		\begin{equation}
				%			\frac{dy}{dt}<0, \qquad 0 < t<\epsilon.
				%		\end{equation}
			%	From inequalities~\eqref{eqs:dy_dt>0<0-ER} it follows that
			%	there exists $0<\tilde t< t_0$ for which $\frac{dy}{dt}(\tilde t)=0$,
			which is a contradiction. Therefore, we conclude that there is no~$t_{0}>0$ for which $\frac{dy}{dt}(t_{0})=0$.  Since we also proved that $\frac{dy}{dt}<0$ for $0 < t \ll 1$, it follows that $\frac{dy}{dt}<0$ for $0 \le t<\infty$.
			
			Since $y(t)$ is monotonically decreasing, $y_{\infty}:=\lim_{t\rightarrow\infty}y$
			can be either a constant or $-\infty$. If we assume by contradiction that $y_{\infty}=-\infty$, 
			then, letting $t\rightarrow\infty$ in~\eqref{eq:y-ER-Bass} gives
			$
			\lim_{t\rightarrow\infty}\frac{dy}{dt}=+\infty,
			$
			which contradicts our earlier result that $\frac{dy}{dt}<0$. Hence, $y_{\infty}$
			is a constant, and so $\lim_{t\rightarrow\infty}\frac{dy}{dt}=0$.
			Letting $t\rightarrow\infty$ in~\eqref{eq:y-ER-Bass} gives
			\[
			0=\lim_{t\rightarrow\infty}\frac{dy}{dt}=-\frac{q}{\lambda} y_{\infty}.
			\]
			Hence, $y_{\infty}=0$. Finally, since $\frac{dy}{dt}<0$, we get that $0<y<1$.
		\end{proof}
	
	Lemma~\ref{lem:f^ER_Bass-monotone} follows from Lemma~\ref{lem:y(t)-properties-ER} and~\eqref{eq:y-ER-Bass}.

	\section{Proof of Lemma~\ref{lem:f^SI-monotone}}
	\label{app:f^ER_SI-monotone}
	
	 we first  show that 
	the auxiliary function $y_{\rm SI}(t)$ is monotonically decreasing in~$t$  towards a positive limit as~$t \to \infty$:
	\begin{lemma}
		\label{lem:y_SI-ER}
		Let  %$q>0$, $p=0$, and $0<I^0<1$, and let 
		$y_{\rm SI}(t)$ be the solution of~\eqref{eq:y-ER-SI}. Then 
		\begin{subequations}
			\begin{equation}
				\label{eq:y_SI(t)-properties-ER}
				\frac{dy_{\rm SI}}{dt}(t)<0, \qquad t>0. 
			\end{equation}
			Let $ y_{\rm SI}^{\infty}:=\lim\limits_{t \to \infty} y_{\rm SI}$. Then $0< y_{\rm SI}^{\infty}<1$. Moreover, 
			$y_{\rm SI}^{\infty}$  
			is the unique root in~$(0, 1)$ of the equation
			\begin{equation}
				\label{eq:y_SI(t)-properties-ER-2}
				y_{\rm SI}^{\infty} = (1-I^0)e^{-\lambda(1-y_{\rm SI}^{\infty})}.
			\end{equation}
		\end{subequations}
		%In particular, $0< y_{\rm SI}^{\infty}<1$.
	\end{lemma}
	\begin{proof}
	We can rewrite equation~\eqref{eq:y-ER-SI} for~$y_{\rm SI}$ as
		$$
		\frac{dy_{\rm SI}}{dt}=\frac{q}{\lambda} h(y), \quad t \ge 0, \qquad  y_{\rm SI}(0)=1, \qquad h(y):=-y+(1-I^0)e^{-\lambda(1-y)}.
		$$
		%where $$. 
		Since 
		$
		h(0) = (1-I^0)e^{-\lambda} >0$ and $h(1) = -I^0<0$,
		$h(y)$ has an odd number of roots in~$(0, 1)$. 
		Assume by contradiction that $h(y)$ has three or more roots. 
		Then by Rolle's theorem, $\frac{dh}{dy}$ has at least two roots 
		in~$(0, 1)$. %\,\footnote{This holds even if the three roots are not distinct, i.e., if  			$h$ has one simple root and one double root.} 
		 Therefore, 
		$\frac{dh^2}{dy^2}$ has at least one root  in~$(0, 1)$.
		This, however, is not possible, since $ \frac{dh^2}{dy^2} =  
		\lambda^2 (1-I^0)e^{-\lambda(1-y)}>0$. Therefore,  $h(y)$~has exactly one root in~$(0, 1)$. 
		Let us denote this root by~$y_{\rm SI}^{\infty}$. 
		Then $y_{\rm SI}^{\infty}$ is the unique solution of~\eqref{eq:y_SI(t)-properties-ER-2} in~$(0, 1)$.
		%	$$
		%	0< y_{\rm SI}^{\infty}<1, \qquad 
		%	y_{\rm SI}^{\infty} = (1-I^0)e^{-\lambda(1-y_{\rm SI}^{\infty})}.
		%	$$

		Since
		$
		\frac{dy_{\rm SI}}{dt}(0) = \frac{q}{\lambda} h(y=1) = -\frac{q}{\lambda}I^0<0,
		$
		and $h(y)$ has a unique root in $[0, 1]$, a standard phase-line analysis shows that
		$$
		\frac{dy_{\rm SI}}{dt}(t)<0, \quad t >0, \qquad  \lim_{t \to \infty} y_{\rm SI}  = y_{\rm SI}^{\infty},
		$$ 
		which concludes the proof.
	\end{proof}

	By~\eqref{eq:f^ER-SI} and~\eqref{eq:y_SI(t)-properties-ER}, 
	$f^{\rm ER}_{\rm SI}(t)$ is monotonically increasing. 
	Letting $t \to \infty$ in~\eqref{eq:f^ER-SI} gives
	$$
	f_{\rm SI}^{\infty}= 1-(1-I^0)e^{-\lambda(1-y_{\rm SI}^{\infty})}= 1- y_{\rm SI}^{\infty}.
	$$	
	where the second equality follows from~\eqref{eq:y_SI(t)-properties-ER-2}.
	Therefore, relation~\eqref{eq:implicit}  follows. 
	Since  $ y_{\rm SI}^{\infty}$ is the unique root in~$(0, 1)$ of~\eqref{eq:y_SI(t)-properties-ER-2}, it follows that $f_{\rm SI}^{\infty}$
	is the unique root in~$(0, 1)$ of~\eqref{eq:implicit}. 
	This concludes the proof of Lemma~\ref{lem:f^SI-monotone}.
	
	 \section{Proof of Corollary~\ref{cor:lim_beta_0_to_0_f_SI^infty=X(lambda)}}
		\label{app:lim_beta_0_to_0_f_SI^infty=X(lambda)}
		
		Letting $I^0 \to 0$ in~\eqref{eq:implicit} gives~\eqref{eq:GC-ER}.
		To finish the proof, we need to show that when $\lambda >1$,  $\lim_{ I^0 \to 0} f_{\rm SI}^{\infty}$ approaches the positive root of~\eqref{eq:GC-ER}, and not the trivial one. 
		To see that, note that for any  $0<I^0 \ll 1$, the giant component has at least one node which is infected at $t=0$. Consequently, as 
		$ t \to \infty$, all the nodes in the giant component become infected.
		Hence, $f_{\rm SI}^{\infty}(I^0) \ge X(\lambda) + I_0(1-X(\lambda)) >X(\lambda)>0$, 
		and hence also in the limit as $I^0 \to 0$.

	\section{Proof of Lemma~\ref{lem:f^d-regular-increasing-in-d}}
	\label{app:f^d-regular-increasing-in-d}
	
	Substituting $t=0$ in~\eqref{eq:S^d-reg} gives
	$\frac{d [S^{ \text{\rm d-reg}}]}{dt}(0) = -(1-I^0)(p+q I^0)$, which 
	is independent of~$d$. 
	Differentiating~\eqref{eq:S^d-reg}  and substituting $t=0$ gives
	$$
	\frac{d^2 [S^{ \text{\rm d-reg}}]}{dt^2}(0)
	= 
	\{\text{terms independent of $d$}\} + \left\{q [S^{ \text{\rm d-reg}}]^2 
	\frac{d}{dt}  \Big(\frac{ e^{pt} [S^{ \text{\rm d-reg}}]}{1-I^0}\Big)^{-\frac{2}{d}}\right\}_{t=0}.
	$$
	Since $q [S^{ \text{\rm d-reg}}]^2(0) =   q (1-I^0)^2$ is positive and independent of~$d$, 
	and since 
	$$
	\begin{aligned} 
		\frac{d}{dt}  \Big(\frac{ e^{pt} [S^{ \text{\rm d-reg}}]}{1-I^0}\Big)^{-\frac{2}{d}} \Big|_{t=0}
		& =-\frac{2}{d} \bigg(\frac{1-I^0}{1-I^0}\bigg)^{-\frac{2}{d}-1} 
		\frac{ p[S^{ \text{\rm d-reg}}](0) +\frac{d [S^{ \text{\rm d-reg}}]}{dt}(0)}{1-I^0} 
		\\ & 
		=-\frac{2}{d} \frac{p(1-I^0)-(1-I^0)(p+q I^0)}{1-I^0}
		=\frac{2}{d}q I^0,
	\end{aligned}
	$$
	we conclude that  $\frac{d^2 [S^{ \text{\rm d-reg}}]}{dt^2}(0)$
	is monotonically decreasing in~$d$.
	Therefore, $[S^{ \text{\rm d-reg}}](t)$ is monotonically decreasing in~$d$
	for $0<t \ll 1$. 
	
	Assume by contradiction that $[S^{ \text{\rm d-reg}}]$ is not monotonically decreasing in~$d$ for all $t>0$. Then there exist $d_1<d_2$ and $t_0>0$ such that 
	$[S^{ \text{\rm $d_1$-reg}}](t_0) \le [S^{ \text{\rm $d_2$-reg}}](t_0)$.
	Therefore, there exists $0<t_1 \le t_0$ such that 
	$[S^{ \text{\rm $d_1$-reg}}](t_1)=[S^{ \text{\rm $d_2$-reg}}](t_1)$
	and 
	$[S^{ \text{\rm $d_1$-reg}}](t)>[S^{ \text{\rm $d_2$-reg}}](t)$ for $0<t <t_1$.
	Hence,
	$\frac{d [S^{ \text{\rm $d_1$-reg}}]}{dt} (t_1) \le \frac{d [S^{ \text{\rm $d_2$-reg}}]}{dt}(t_1)$.
	Since, however, $[S^{ \text{\rm $d_1$-reg}}](t_1)=[S^{ \text{\rm $d_2$-reg}}](t_1)<(1-I^0)e^{-pt}$, it follows from~\eqref{eq:S^d-reg} that $\frac{d [S^{ \text{\rm $d_1$-reg}}]}{dt} (t_1) > \frac{d [S^{ \text{\rm $d_2$-reg}}]}{dt}(t_1)$, which is a contradiction. Therefore, we conclude that
	$[S^{ \text{\rm d-reg}}]$ is monotonically decreasing in~$d$ for all times,
	and so the result follows from~\eqref{eq:f^d-regular=1-[S^d-reg]}. 
	
%	 Let $t>0$.  From~\eqref{eqs:multisub-1} and the fact that 
%	   $[S_m] > [S_{k_i,m}]$ we have that $[S_m]<(1-I^0) e^{-pt}$ . 
%	Therefore, $[S] =\frac1M \sum_{m=1}^M [S_m]<(1-I^0) e^{-pt}$,  and so 
%		% 	
%		% 	Let $z(t):=\frac{e^{pt} }{1-I^0}[S]$, where $[S]$ is the solution of~\eqref{eq:S^d-reg}.   Then
%		%	\begin{equation}
%			%		\label{eq:S^d-regular-z}
%			%	\frac{dz}{dt}= -q z \left( 1- (1-I^0) e^{-pt} z^{-\frac{2}{d}} z \right), \qquad z(0 )= 1.
%			%\end{equation}
%			%Since $z<1$, 
%			the right-hand side of~\eqref{eq:S^d-reg} is monotonically decreasing in~$d$. Hence, %by Lemma~\ref{lem:dif_ineq}, 
%			$[S^{ \text{\rm d-reg}}]$ is monotonically decreasing in~$d$,
%			and so the result follows from~\eqref{eq:f^d-regular=1-[S^d-reg]}.
%%		
%	Let $z(t):=\frac{[S]}{e^{-pt}(1-I^0)}$, where $[S]$ is the solution of~\eqref{eq:S^d-reg}.   Then
%	\begin{equation}
%		\label{eq:S^d-regular-z}
%		\frac{dz}{dt}= -q z \left( 1- (1-I^0) e^{-pt} z^{-\frac{2}{d}} z \right), \qquad z(0 )= 1.
%	\end{equation}
%	Since $z<1$, the right-hand side of~\eqref{eq:S^d-regular-z}  is monotonically decreasing in~$d$. Therefore, %by Lemma~\ref{lem:dif_ineq}, 
%	$z(t)$ is 
%	monotonically decreasing in~$d$. Hence, so is $[S^{ \text{\rm d-reg}}]$.
%	
	\section{Proof of Lemma~\ref{lem:bounds-f^regular}}
	\label{app:bounds-f^regular}
	
	When $d=2$,  equation~\eqref{eq:S^d-reg} reads 
	\begin{equation*}
		\frac{d [S]}{dt} = -[S] \left( p+q \left(1-(1-I^0) e^{-pt}  \right) \right), \qquad [S](0 )= 1-I^0.
	\end{equation*}
	This equation coincides with that for~$[S^{\rm 1D}] = 1-f^{\rm 1D}$, see~\cite{OR-10,Bass-monotone-convergence-23}.
	Therefore, we have~\eqref{eq:f^2-regular_equiv_f^1D}.
	Letting $d \to \infty$ in equation~\eqref{eq:S^d-reg}, we get 
	\begin{equation*}
		\frac{d [S]}{dt} = -[S] \left( p+q \left(1-  [S] \right) \right), \qquad [S](0 )= 1-I^0.
	\end{equation*}
	This is exactly the ODE for $[S^{\rm compart}] = 1-f^{\rm compart}$,
	see~\eqref{eq:f^complete_infty-ODE}. 
	Finally, since $f^{ \text{\rm d-reg}}$ is monotonically increasing in~$d$
	(Lemma~\ref{lem:f^d-regular-increasing-in-d}),    inequalities~\eqref{eq:bounds-f^regular} follow.

	 \section{Proof of Lemma~\ref{lem:E_G[c_L^d]-bound}}
	\label{app:E_G[c_L^d]-bound}
	
		Let $j$ belong to a cycle of length~$L$. Denote the cycle nodes by 
	$\{j,m_1, \dots, m_{L-1}\}$. There are $M-1$ possibilities for~$m_1$,
	$M-2$ possibilities for~$m_2$, up to $M-(L-1)$ possibilities for~$m_{L-1}$.
	Therefore, there are $\prod_{i=1}^{L-1}(M-i)$ such cycles.
	We should, however, divide this product by~2, because each cycle is counted both clockwise and counter-clockwise. 
	
	Since $j$ has degree $d$, the probability in~${\cal G}(M,\alpha)$ that $j$ is connected to~$m_1$ and to $m_{L-1}$ is $\frac{d}{M-1} \frac{d-1}{M-2}$. For $i=1, \dots, L-2$, the probability that $m_i$ is connected to~$m_{i+1}$ is~$\alpha$.
	Therefore, the probability that the cycle $\{j,m_1, \dots, m_{L-1}\}$ exists in~${\cal G}(M,\alpha)$ is
	$\frac{d}{M-1} \frac{d-1}{M-2} \alpha^{L-2}$. 
	
	Combining the above and using~$\alpha  =  \frac{\lambda}{M}$, we have 	  
	$$
	C_L^d  =
	\frac{d}{M-1} \frac{d-1}{M-2} \alpha^{L-2} \frac12 \prod_{i=1}^{L-1}(M-i) 
	= \frac{\kappa_{L,d}^{\rm ER}}{M} \prod_{i=3}^{L-1}\frac{M-i}{M}. 
	% =  d(d-1)(M-3) \dots (M-(L-1)) \frac1{2} \alpha^{L-2},
	$$
	Therefore, relations~\eqref{eq:E_G[c_L^d]-bound} and~\eqref{eq:E_G[c_L^d]} follow.

\section{Proof of Lemma~\ref{lem:upper-bound-c_L^d-regular}}
\label{app:upper-bound-c_L^d-regular}

	In a complete graph with $M$~nodes,
the number of different cycles of size~$L$  is
$
\frac{M(M-1) \cdots (M-L+1)}{2L}.
$
This is because there are $M$~possibilities for the ``first'' node, $M-1$~for the second,  and up to 
$M-(L-1)$ for the ``last'' node. In this count, however, each cycle is counted $L$~times by translation and 2~times by reflection. Therefore, 
we need to divide the numerator by~$2L$. 

Let $L \in \{3, \dots, M\}$.
By~\cite[Corollary~2.15]{Ben-Shimon-11}, there exists $\beta_d>0$, such that the probability that any given cycle of length~$L$ 
exists in ${\cal G}_{M}^{\text{\rm d-reg}}$
is bounded from above by $\left(\frac{\beta_d}{M}\right)^L$. 
Therefore, the expected number of different cycles of size~$L$ in a $d$-regular graph with $M$~nodes 
has the upper  bound
$$
\mathbb{E}_{{\cal G}_{M}^{\text{\rm d-reg}}}
\big[\text{\# of cycles of length $L$ in $ G^{\text{\rm d-reg}}$}  \big]
< \frac{M(M-1) \cdots (M-L+1)}{2L}
\left(\frac{\beta_d}{M}\right)^L 
< \frac{(\beta_d)^L}{2L} . 
$$
The probability that $j$ belongs to a cycle of size~$L$ is~$\frac{L}{M}$. Therefore, 
$$
C_L^{ \text{\rm d-reg}} = 
\frac{L}M \mathbb{E}_{{\cal G}_{M}^{\text{\rm d-reg}}}\big[\text{\# of cycles of length $L$ in $ G^{\text{\rm d-reg}}$}  \big],
$$
and so the result follows.

%\section{Proof of Lemma~\ref{lem:c^d-regular>3(d/2-1)}}
%\label{app:c^d-regular>3(d/2-1)}
%
%We first prove an auxiliary lemma:
%\begin{lemma}
%	\label{lem:Kriv}
%	Let $G$ be a graph with $M$ nodes and $|E(G)|$~edges. Then 
%	%for any $d$-regular graph on $M$~vertices, 
%	$$
%	\text{\rm \# of cycles in $G$} \, \ge |E(G)|-(M-1).
%	$$
%\end{lemma}
%\begin{proof}
%Denote by~$F$
%the {\em spanning forest} of~$ G$.\,\footnote{A {\em 
%		spanning forest}~$F$ of an undirected graph~$G$ is a subgraph that includes all the vertices of~$G$ and a subset of its edges, such that for each connected component of~$G$,  all the nodes in~$F$ are connected and there are no cycles.} Then $F$ contains at most 
%$M-1$~edges.
%Every edge of $G-E(F)$ closes a cycle in~$G$ with~$F$, and these cycles are different for different edges. Hence, there are at least $|E(G)|-|E(F)|\ge |E(G)|-(M-1)$ different cycles.
%\end{proof}
%
%
%\textbf{Proof of Lemma~\ref{lem:c^d-regular>3(d/2-1)}.~}
%Since $|E(G^{\rm regular})| = \frac{dM}2$, it follows from Lemma~\ref{lem:Kriv}
%that 
%
%$$
%\text{\rm \# of cycles in $G^{\rm regular}$} >
%\frac{dM}{2}-(M-1) > M \left(\frac{d}{2}-1 \right). 
%$$
%Each cycle adds $L \ge 3$ nodes that lie on a cycle. 
%Therefore, 
%$$
%c^{ \text{\rm d-reg}}( G^{\rm regular}) \, \ge  \left(\text{\rm \# of cycles in $G^{\rm regular}$}\right) \frac3{M}
%\ge  3 \left(\frac{d}{2}-1 \right).
%$$ 

%\end{APPENDICES}

	\bibliographystyle{plain}
	\bibliography{../../diffusion}

\end{document}